\documentclass[a4paper,10pt]{scrartcl} 

\usepackage{graphicx}
\usepackage{overpic}         
\usepackage{amsfonts}
\usepackage{amsmath}         
\usepackage{amssymb}	     
\usepackage[usenames,dvipsnames]{color} 
\usepackage{algorithm}       
\usepackage{algpseudocode}   
\usepackage{url}             
\usepackage{arydshln}        
\usepackage{multirow}        

\usepackage{bm} 

\usepackage{etoolbox}    
\let\bbordermatrix\bordermatrix
\patchcmd{\bbordermatrix}{8.75}{4.75}{}{}
\patchcmd{\bbordermatrix}{\left(}{\left[}{}{}
\patchcmd{\bbordermatrix}{\right)}{\right]}{}{}

\usepackage{tcolorbox}
\usepackage{colortbl}
\usepackage{geometry}

\usepackage{listings,color}
\definecolor{mygray}{gray}{0.5}
\definecolor{myred}{rgb}{0.8,0,0}
\definecolor{myblue}{rgb}{0,0,0.6}
\definecolor{mybrown}{rgb}{0.1,0.5,0.5}

\usepackage{caption} 
\usepackage{subcaption} 

\usepackage{paralist} 

\newcommand{\lambdaplus}{\mathop{\lambda}^{\textnormal{\tiny{+}}}{}}
\newcommand{\lambdaminus}{\mathop{\lambda}^{\textnormal{\tiny{--}}}{}}
\newcommand{\funcplus}{\mathop{u}^{\textnormal{\tiny{+}}}{}}
\newcommand{\funcminus}{\mathop{u}^{\textnormal{\tiny{--}}}{}}
\newcommand{\deltaRef}{{\delta}_j^{(h)}}
\newcommand{\deltaHAMLSBEM}{\widehat{\delta}_j^{\, (h)}}

\newcommand{\diag}{\mathrm{diag}}

\newcommand{\Nidx}{^{\textnormal{\tiny{$(N)$}}}}
\newcommand{\hidx}{^{(h)}}

\newcommand{\Nh}{N_h}

\newcommand{\Ocal}{\mathcal{O}}
\newcommand{\setminusNull}{\setminus \{0\}}
\newcommand{\Id}{\mathrm{Id}}
\newcommand{\dx}{\, \, \mathrm{d}x}
\newcommand{\dy}{\, \, \mathrm{d}y}

\newcommand{\nev}{n_{\textnormal{es}}}
\newcommand{\divop}{\mathrm{div}}
\newcommand{\R}{\mathbb{R}}
\newcommand{\N}{\mathbb{N}}
\newcommand{\C}{\mathbb{C}}

\newcommand{\LDL}{LDL$^{\textnormal{T}}$}
\newcommand{\LU}{LU}

\newcommand{\Ssum}{Z}
\newcommand{\ksum}{\bar{k}}

\newcommand{\HM}{\ensuremath{\mathcal{H}}}

\newcommand{\Rk}{\ensuremath{\mathrm{R}(k)}}
\newcommand{\task}{\ensuremath{\mathbb{T}}}
\newcommand{\ctask}{\ensuremath{\mathbb{CT}}}

\newcommand{\MaxDofSubdomain}{N^\textnormal{\tiny{AMLS}}_{\min}}

\newcommand{\rmax}{l_\textnormal{\tiny{rec}}}

\newcommand{\Ph}{\mathcal{P}}

\newcommand{\LTwo}{L^{2}(\Omega)}
\newcommand{\LSup}{L^{\infty}(\Omega)}

\newcommand{\Lightblue}{\textcolor{blue!80!white}}

\DeclareOldFontCommand{\bf}{\normalfont\bfseries}{\mathbf}
\newtheorem{theorem}{Theorem}[section]

\newtheorem{remark}[theorem]{Remark}
\newtheorem{definition}[theorem]{Definition}

\newenvironment{proofsquare}{\noindent {\bf Proof:}}{ {\hfill{\vrule height 6pt width 6pt depth 0pt}}\newline}

\evensidemargin=\oddsidemargin
\newdimen\dummy
\dummy=\oddsidemargin
\begin{document}

\title{
Solving an Integral Equation Eigenvalue Problem
via a New Domain Decomposition Method and Hierarchical Matrices}
\subtitle{{\textnormal{
Automated Multi-Level Substructuring for Dense Eigenvalue Problems
}}}

\author{
  Peter Gerds\thanks{
    Institut f\"ur Geometrie und Praktische Mathematik, 
    RWTH Aachen University, Templergraben 55, 52056 Aachen, Germany. 
    Email: {\texttt{gerds@igpm.rwth-aachen.de}}.
    }
}

\date{}

\maketitle

\begin{abstract}
In this paper the author introduces a new domain decomposition method for the solution of discretised
integral equation eigenvalue problems. 
The new domain decomposition method is motivated by the so-called \emph{automated multi-level substructuring} (short 
AMLS) method. 
The AMLS method is a domain decomposition method for the solution of elliptic PDE eigenvalue problems
which has been shown to be very efficient especially when a large number of eigenpairs is sought.
In general the AMLS method is only applicable 
to these kind of continuous eigenvalue problems 
where the corresponding discretisation leads to an algebraic eigenvalue problem of the form 
$Kx=\lambda Mx$ where $K,M \in \R^{N\times N}$ are symmetric sparse matrices.
However, the discretisation of an integral equation eigenvalue problem leads to a discrete problem 
where the matrix $K$ is typically dense, since a non-local integral operator is involved in the equation. 
In the new method, which is introduced in this paper,
the domain decomposition technique of classical AMLS is generalised to eigenvalue problems $Kx=\lambda Mx$ where $K,M$ 
are symmetric and possibly dense, and which is therefor applicable for the solution of integral equation eigenvalue 
problems. To make out of this an efficient eigensolver, the new domain decomposition method is combined with a 
recursive approach and with the concept of {hierarchical matrices}. \\



\end{abstract}

\section{Introduction}
\label{section: Introduction}

Integral equation eigenvalue problems are typically solved using the finite element discretisation
where the arising discrete eigenvalue problem is solved by an iterative algebraic eigensolver or a direct method 
(e.g., QR-algorithm). 
Domain decomposition techniques, which are available for example for the 
solution of elliptic PDE eigenvalue problems (see, e.g., \cite{Bennighof_Lehoucq,Seshu}), however, 
are not known (to the best of the author's knowledge)
for integral equation eigenvalue problems.
In this paper the very first domain decomposition method for the solution of integral equation 
eigenvalue problems is introduced. This new domain decomposition method is motivated 
by the so-called \emph{automated multi-level substructuring} (short AMLS) method. 
The AMLS method was mainly developed by Bennighof and co-authors 
\cite{First_Bennighof,Bennighof_Lehoucq,Kaplan} for the solution of elliptic PDE eigenvalue problems,
and it has been shown that AMLS is very efficient especially when a large number of eigenpairs is sought
\cite{Bennighof_Extending_Frequency_Response,Kaplan,Kropp_Heiserer}. In general the AMLS method is only applicable 
to these kind of continuous eigenvalue problems where the corresponding discretisation leads to an algebraic eigenvalue 
problem of the form $Kx=\lambda Mx$ where the matrices $K,M \in \R^{N\times N}$ are 
symmetric sparse.
However, the discretisation of integral equation eigenvalue problems leads to a discretised problem
$Kx=\lambda Mx$
where the matrix $K$ is typically dense, because a non-local integral operator is involved in the 
equation.
In the new method, which is introduced in the following 
and which is called \emph{dense AMLS method}, 
the domain decomposition technique of classical AMLS is generalised to eigenvalue problems $Kx=\lambda Mx$ where $K,M$
are symmetric and possibly dense, and which is therefore applicable for the solution of integral equation eigenvalue 
problems. To improve the efficiency of this new domain decomposition technique, the dense AMLS method
is combined with a recursive approach 
and with the concept of the \emph{hierarchical matrices} (short \HM-matrices).

The remainder of this paper is organised as follows:
In Section \ref{section: Problem Description}
the integral equation eigenvalue problem and the underlying problem setting is introduced.
After this, in Section \ref{section: AMLS classical}, the basic steps of the classical AMLS method are summarised.
In Section \ref{section: AMLS dense} the new dense AMLS method is introduced
and in Section \ref{section: Numerical Results} numerical results are presented.
Finally, in Section \ref{section: Recursive Approach and Use of Hierarchical Matrices}
a recursive version of dense AMLS is introduced which is combined with the concept of \HM-matrices. 
A summary and conclusion of this work is given in Section \ref{section: Conclusion}.

\section{Problem Description}
\label{section: Problem Description}

In this paper we consider the continuous eigenvalue problem
\begin{equation}
(\mathop{A} u)(x) = \lambda u(x) 
\quad \textnormal{for all } x \in \Omega
\label{introduction: continuous EVP classic}
\end{equation}
where $\Omega \subset \R^d$ is an open bounded domain ($d \in \N$),
$\mathop{A}$ is a linear Fredholm integral operator of the form 
\begin{equation}
(\mathop{A} u)(x) 
:= 
\int_{\Omega}
k(x,y) u(y) \dy
\quad \textnormal{with } x \in \Omega
\label{introduction: integral operator}
\end{equation}
with kernel function $k: \Omega \times \Omega \to \R$, 
and where in \eqref{introduction: continuous EVP classic} a suitable eigenfunction $u: \Omega \to \C$ and the 
associated eigenvalue 
$\lambda \in \C$ are sought.
Throughout the paper the following definitions are used 
for the kernel function $k$:
\begin{itemize}
\item 
$k$ is called \emph{continuous} iff $k \in C^0(\overline{\Omega} \times \overline{\Omega})$.

\item 
$k$ is called
\emph{weak singular} iff
$k \in C^0( \{ (x,y) \in \overline{\Omega} \times \overline{\Omega} : x \neq y \} )$
and there exist constants $C>0$ and  $\alpha \in [0,d)$ such that
\begin{equation*}
|k(x,y)| \leq C |x-y|^{-\alpha}
\qquad \textnormal{for all }
x,y \in \overline{\Omega} \textnormal{ with } x \neq y.
\end{equation*}

\item
$k$ is called
\emph{symmetric} iff $k(x,y)=k(y,x)$ for all 
$x,y \in \overline{\Omega} \textnormal{ with } x \neq y$.

\end{itemize}
From the theory of linear integral equations (see, e.g., \cite{Hackbusch_integralEQUATIONS}) the following result is 
known:
\begin{theorem}[Existence of Eigensolutions]
\label{introduction: theorem: existence eigensolution}
Let the kernel function $k$ in \eqref{introduction: integral operator}
be continuous or weak singular with $\alpha \in [0,d/2)$, let $k$ be symmetric and assume that $k \neq 0$. Then the 
integral 
operator 
$\mathop{A}: \LTwo \to \LTwo$ in \eqref{introduction: integral operator}
is selfadjoint and compact, and hence
eigenvalue problem \eqref{introduction: continuous EVP classic} possesses a countable 
family of eigensolutions
\begin{equation}
\bigl (
\lambda_j,u_j
\bigr )_{j=1}^{\infty} \in \mathbb{R} \times \LTwo \setminusNull 
\label{introduction: continuous eigensolution}
\end{equation}
where all eigenvalues $\lambda_j$ are ordered 
with respect to their multiplicity%
\footnote{This means that the eigenvalues in \eqref{introduction: continuous eigensolution} are repeated 
according to the dimension of the corresponding eigenspace}
and absolute value such that $|\lambda_1| \geq |\lambda_2| \geq \ldots \geq 0$.
In particular, it holds
\begin{itemize}
\item[i)]
All eigenvalues $\lambda_j$ are real and 
we have $\lambda_j \xrightarrow{\; \; j \rightarrow \infty \; \;} 0$.

\item[ii)]
The eigenspace $E(\lambda_j) \subset \LTwo$ of the eigenvalue $\lambda_j$, which is defined by
\begin{equation}
E(\lambda_j)
:=
\mathrm{span}
\Bigl \{
u \in \LTwo
\; : \;
\mathop{A} u= \lambda_j u
\Bigr \},
\label{introduction: continuous eigenspace from existence theorem}
\end{equation}
is finite-dimensional for $\lambda_j \neq 0$.

\item[iii)]
If it holds $\lambda_j \neq \lambda_k$ then the corresponding eigenfunctions
$u_j$ and $u_k$ are orthogonal with respect to the $\LTwo$-inner product, i.e.,
it holds $(u_j, u_k)_0 := \int_{\Omega} u_j(x) u_k(x) \dx = 0$.

\item[iv)]
The eigenfunctions $\bigl (u_j \bigr )_{j=1}^{\infty}$ form a basis of the Hilbert space $\LTwo$ 
and without loss of generality it can be assumed that all eigenfunctions are orthonormal with respect to 
$(\cdot,\cdot)_0$.

\end{itemize}
\end{theorem}

Multiplying \eqref{introduction: continuous EVP classic} by $v \in \LTwo$ and integrating over $\Omega$ we obtain 
the equivalent variational eigenvalue problem 
\begin{equation}
\begin{cases}
\; \textnormal{find } \, (\lambda,u) \in  \R \times \LTwo \setminusNull \, \textnormal{ such that}
\\
\; a(u,v) \, = \, \lambda \, (u,v)_0 
\qquad \forall \:
v \in \LTwo
\end{cases}
\label{introduction: continuous EVP variational}
\end{equation}
with the symmetric bilinear form $a(u,v):=\int_{\Omega}\int_{\Omega}  v(x) k(x,y) u(y) \dy \dx$.
We approximate solutions of the continuous eigenvalue problem 
\eqref{introduction: continuous EVP classic}
and \eqref{introduction: continuous EVP variational}
by discretisation:
Using a conforming finite element space $V_h \subset \LTwo$
with dimension $N_h$ and nodal basis
$
\bigl (
\varphi_i\hidx
\bigr )_{i=1}^{N_h}
$
the eigenvalue problem \eqref{introduction: continuous EVP variational} is discretised by
\begin{equation}
\begin{cases} 
\; \textnormal{find }  (\lambda\hidx, x\hidx  ) \in  \R \times \R^{\Nh} \setminus \{0\} \textnormal{ with} 
\rule[-3mm]{0mm}{3mm}
\\
\; K\hidx \, x\hidx \, = \, \lambda\hidx \, M\hidx \, x\hidx
\end{cases}
\label{introduction: discrete EVP}
\end{equation}
where the stiffness and mass matrix
\begin{equation}
K\hidx
:=  
\Bigl ( a  ( \varphi_j\hidx, \varphi_i\hidx )  \Bigr )_{i,j=1}^{\Nh} 
\hspace*{-1mm}
\in \R^{\Nh \times \Nh}
\quad \textnormal{and} \quad 
M\hidx
:= 
\Bigl (   (  \varphi_j\hidx, \varphi_i\hidx   )_0 \Bigr )_{i,j=1}^{\Nh}
\hspace*{-1mm}
\in \R^{\Nh \times \Nh}
\label{introduction: stiffness and mass  matrix}
\end{equation}
are both symmetric. In contrast to the 
sparse and positive definite mass matrix $M$, the stiffness matrix $K$ is in general dense 
(i.e., nearly all entries of $K$ are nonzero) and in general not positive definite. 

From the approximation theory of integral equation eigenvalue problems follows that the discrete eigensolutions
$(\lambda\hidx,u\hidx):=(\lambda\hidx,\Ph x\hidx)$
are approximating 
the continuous eigensolutions $(\lambda,u)$ of 
\eqref{introduction: continuous EVP classic}, 
where $u\hidx$ is associated to the discrete problem \eqref{introduction: discrete EVP} via 
\begin{equation*}
\Ph: \R^{\Nh} \rightarrow V_h \subset \LTwo
\quad \textnormal{with} \quad
x\hidx \mapsto \sum_{i=1}^{\Nh} x_i\hidx \, \varphi_i\hidx 
\,.
\label{introduction: prolongation operator FEM}
\end{equation*}
For the precise formulation of the corresponding approximation result an
additional notation of the continuous and discrete eigensolutions is introduced which distinguishes between positive 
and 
negative eigensolutions:
In the following eigensolutions of the continuous and discrete problem associated to the positive eigenvalues are 
denoted by
\begin{alignat*}{2}
\bigl (
\lambdaplus_j,\funcplus_j
\bigr )
&
\in \mathbb{R} \times \LTwo \setminusNull 
& & \qquad \textnormal{where} \quad
\lambdaplus_1 \geq \lambdaplus_{2} \geq \ldots \geq 0,
\\
\bigl (
\, \lambdaplus_j\hidx\,  , \funcplus_j\hidx \,
\bigr )
&
\in \R \times V_h \setminusNull
& & \qquad \textnormal{where} \quad
\lambdaplus_1\hidx \geq \lambdaplus_{2}\hidx \geq \ldots \geq 0
\end{alignat*}
and eigensolutions associated to the negative eigenvalues by
\begin{alignat*}{2}
\bigl (
\lambdaminus_j,\funcminus_j
\bigr )
&
\in \mathbb{R} \times \LTwo \setminusNull 
& & \qquad \textnormal{where} \quad
\lambdaminus_1 \leq \lambdaminus_{2} \leq \ldots \leq 0,
\\
\bigl (
\, \lambdaminus_j\hidx\,  , \funcminus_j\hidx \,
\bigr )
&
\in \R \times V_h \setminusNull
& & \qquad \textnormal{where} \quad
\lambdaminus_1\hidx \leq \lambdaminus_{2}\hidx \leq \ldots \leq 0.
\end{alignat*}
\begin{theorem}[Qualitative Convergency Results] 
\label{introduction: theorem: Ritz-Galerkin Approximation Qualitative Results}
Let the assumptions of Theorem \ref{introduction: theorem: existence eigensolution} be valid. Consider eigenvalue 
problem \eqref{introduction: continuous EVP classic} and its finite element discretisation
\eqref{introduction: discrete EVP} where the eigensolutions are indexed as above, and where 
the finite element space $V_h \subset \LTwo$ fulfils approximation property
\begin{equation}
\lim_{h\rightarrow 0} 
\;
\inf_{v\hidx \in V_h} \|u-v\hidx\|_0
\;
=
0
\qquad 
\textnormal{for all } u \in \LTwo
\label{introduction: sequence of Galerkin spaces}
\end{equation} 
with $\|\cdot\|_0:=(\cdot,\cdot)_0^{1/2}$.
Then the discrete eigenvalues are approximating the continuous ones, i.e., it holds
\begin{equation*}
\lambdaplus_j\hidx
\xrightarrow{\; \; h \rightarrow 0 \; \;} 
\lambdaplus_j
\quad \textnormal{and} \quad 
\lambdaminus_j\hidx
\xrightarrow{\; \; h \rightarrow 0 \; \;} 
\lambdaminus_j
\qquad
\textnormal{for } j \in \N;
\end{equation*}
and for the associated discrete eigenfunctions 
$\displaystyle \funcplus_j\hidx$ and $\displaystyle \funcminus_j\hidx$ 
[assuming that $\displaystyle \|\funcplus_j\hidx\|_0=1$ and $\displaystyle \|\funcminus_j\hidx\|_0=1$] 
exist subsequences
which converge in $\LTwo$ to an eigenfunction $\displaystyle \funcplus \in E(\lambdaplus_j)$
and to an eigenfunction $\displaystyle \funcminus \in E(\lambdaminus_j)$
when $h \to 0$.
\end{theorem}
\begin{proofsquare}
The result for the convergence of the eigenvalues follows directly from the well-known Minimum-Maximum principle
\begin{align*}
\lambdaplus_j
=
\hspace{-1mm}
\mathop{\max}_{
\substack{H_j \subset  \LTwo,  \\  \dim H_j = j}
}
\;
\mathop{\min}_{
\substack{u \in H_j \setminusNull}
}
\dfrac{a(u,u)}{(u,u)_0}
\quad \textnormal{and} \quad
\lambdaplus_j\hidx
=
\hspace{-1mm}
\mathop{\max}_{
\substack{H_j \subset  V_h, \\  \dim H_j = j}
}
\;
\mathop{\min}_{
\substack{u \in H_j \setminusNull}
}
\dfrac{a(u,u)}{(u,u)_0}
\quad \textnormal{for } j=1,2,\ldots,
\\
\lambdaminus_j
=
\hspace{-1mm}
\mathop{\min}_{
\substack{H_j \subset  \LTwo,  \\  \dim H_j = j}
}
\;
\mathop{\max}_{
\substack{u \in H_j \setminusNull}
}
\dfrac{a(u,u)}{(u,u)_0}
\quad \textnormal{and} \quad
\lambdaminus_j\hidx
=
\hspace{-1mm}
\mathop{\min}_{
\substack{H_j \subset  V_h, \\  \dim H_j = j}
}
\;
\mathop{\max}_{
\substack{u \in H_j \setminusNull}
}
\dfrac{a(u,u)}{(u,u)_0}
\quad \textnormal{for } j=1,2,\ldots
\end{align*}
related to compact selfadjoint operators, 
from the continuity of the Rayleigh quotient 
$R(u):=
a(u,u)/(u,u)_0
$, and from the approximation property \eqref{introduction: sequence of Galerkin spaces} of $V_h$.
The result for the convergence of the eigenfunctions is a combination of 
\cite[Satz 4.8.15]{Hackbusch_integralgleichungen} and the convergence result of the eigenvalue approximation described 
above.
\end{proofsquare}

When eigensolutions $(\lambda_j,u_j)$ of integral equation eigenvalue problems are approximated by the finite element 
method the following point has to be noted: Only the by magnitude largest 
eigenvalues $\lambda_j$ and their corresponding eigenfunctions $u_j$ can be well approximated by the finite element 
space $V_h$ because the approximation error increases with by magnitude decreasing eigenvalue
(see, e.g., numerical results from Section \ref{section: Numerical Results}). 

Hence, in the following we are interested only in computing a portion of the eigenpairs of 
the discrete problem \eqref{introduction: discrete EVP}, e.g., the by magnitude largest 
\begin{equation}
\nev = C  N_h^{1/3} \in \N
\quad
\textnormal{ or }
\quad
\nev = C  N_h^{1/2} \in \N
\label{introduction: number of sought eigensolutions}
\end{equation}
eigenpairs (with some $C>0$).

Depending on the number of sought eigenpairs, different approaches are better suited for 
the solution of the algebraic eigenvalue problem \eqref{introduction: discrete EVP}:
If the number of sought eigenpairs $\nev$ is rather small, e.g., if $\nev = 5$,
an iterative algebraic eigensolver such as the Lanczos method \cite{Templates_EWP,Shift_Invert_Lanczos}
is a good choice for the solution of \eqref{introduction: discrete EVP}. 
If the number of sought eigenpairs approaches $N_h$, it is
advisable to use instead a cubic scaling direct method such as the QR algorithm \cite{QR}.

However, since we are interested in a large number of eigenpairs
where $\nev$ is, e.g., of the size
\eqref{introduction: number of sought eigensolutions},
neither an iterative eigensolver nor a direct method might be a good choice for the solution
of the dense problem \eqref{introduction: discrete EVP}.
Using domain decomposition techniques might be a new approach for the efficient solution of  
\eqref{introduction: discrete EVP}.
The domain decomposition techniques, which is applied in the AMLS method
for the solution of elliptic PDE eigenvalue value problems, has been shown to be 
very efficient, especially when the number of sought eigenpairs is, e.g., of the size
\eqref{introduction: number of sought eigensolutions} [see, e.g., \cite{Diss_Gerds,Gerds}]. 
The efficiency of the AMLS method for the solution of elliptic PDE eigenvalue problems
motivates to transfer the domain decomposition technique of AMLS to integral equation eigenvalue problems. 
For doing this the domain decomposition technique of classical AMLS is briefly described in the following section.

\section{The Classical AMLS Method for Sparse Eigenvalue Problems}
\label{section: AMLS classical} 

This section summarises the basic steps of the classical AMLS method for the 
solution of an elliptic PDE eigenvalue problem. The following discussion is restricted to a purely 
algebraic setting and only to a single-level version of the AMLS method. To see a multi-level version of AMLS and to 
see how AMLS is motivated in the context of a continuous setting, it is referred, e.g., 
to \cite{Bennighof_Lehoucq,Gerds}.

The initial point of this section is the 
elliptic PDE eigenvalue problem 
\begin{equation}
\begin{cases}
\; {L} u = \lambda u \quad \textnormal{in } \Omega,
\\
\; \phantom{L} u = 0 \phantom{\lambda} \quad \textnormal{on } \partial\Omega
\end{cases}
\label{AMLS classic: continuous EVP classic}
\end{equation}
where $\Omega \subset \R^d$ is a bounded Lipschitz domain and 
$L$ is a uniformly elliptic second order partial differential operator in divergency form 
$Lu
=
-\divop
\bigl (
A \nabla u
\bigr )
+
c u 
$ with
$A \in (\LSup)^{d \times d}$, $c \in \LSup$ and $c\geq0$.
We approximate solutions of the continuous eigenvalue problem 
\eqref{AMLS classic: continuous EVP classic}
by discretisation.
Using a conforming $N$-dimensional finite element space $V_h$
the eigenvalue problem \eqref{AMLS classic: continuous EVP classic}, respectively its weak formulation, 
is discretised by
\begin{equation}
\begin{cases} 
\; \textnormal{find } \left (\lambda, x \right ) \in  \R \times \R^N \setminusNull  \textnormal{ with}
\\
\; K \, x \, = \, \lambda \, M \, x
\end{cases}
\label{AMLS classic: discrete EVP}
\end{equation}
where the stiffness matrix $K \in \R^{N\times N}$ and the mass matrix 
$M \in \R^{N\times N}$
are sparse and symmetric positive definite, and 
where the 
eigenpairs 
\begin{equation}
\bigl (
\lambda_j,x_j
\bigr )_{j=1}^{N} \in \R_{>0} \times \R^N \setminusNull
\qquad \textnormal{with} \quad
\lambda_j \leq \lambda_{j+1}
\label{AMLS classic: discrete eigenpairs}
\end{equation}
of problem \eqref{AMLS classic: discrete EVP} have all positive eigenvalues.
For reasons of convenience
the upper index of $\lambda\hidx, x\hidx, K\hidx$ and $M\hidx$,
which is indicating the mesh width $h$ of the underlying finite element discretisation,
is left out in \eqref{AMLS classic: discrete EVP}.

The discrete eigenpairs $(\lambda_j,x_j)$ 
of problem \eqref{AMLS classic: discrete EVP}
provide approximations of the eigensolutions of problem 
\eqref{AMLS classic: continuous EVP classic}. 
However, it has to be noted that
only the discrete eigenpairs associated to the smallest 
eigenvalues provide good approximations of continuous eigensolutions since
in the context of elliptic PDE eigenvalue problems
the error of the finite element approximation increases with the size of the eigenvalue 
(see, e.g., \cite{Sauter_Experimental_Study,Sauter_Error_Estimates}).
Hence, we are interested only in the eigenpairs
\eqref{AMLS classic: discrete eigenpairs} which are associated to the smallest eigenvalues,
where the number of sought eigenpairs $\nev$ is, e.g., of the size
\eqref{introduction: number of sought eigensolutions}.
\\

\begin{figure}[t]
\begin{subfigure}[t]{0.47\textwidth}
\centering
      \begin{overpic}[width=6.0cm]{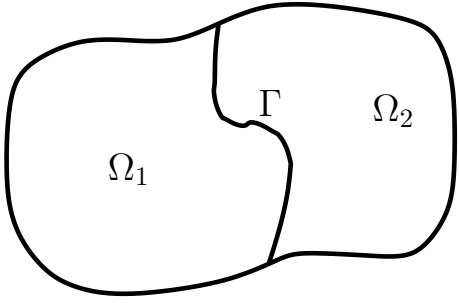}
	\put(23,26){\Large {$\Omega_1$}}
	\put(81,39){\Large {$\Omega_2$}}
	\put(56,40){\Large {$\Gamma$}}
      \end{overpic}
\caption{
Substructuring of the domain $\Omega$ into the non-overlapping subdomains $\Omega_1$ and $\Omega_2$,
where $\overline{\Omega} = \overline{\Omega}_1 \cup \overline{\Omega}_2$,
$\Omega_1 \cap \Omega_2 = \emptyset$
and $\Gamma := \overline{\Omega}_1 \cap \overline{\Omega}_2$.}
\label{AMLS classic: subfigure: Omega substructuring}
\end{subfigure}
\hfill
\begin{subfigure}[t]{0.47\textwidth}
\centering
      \begin{overpic}[width=6.0cm]{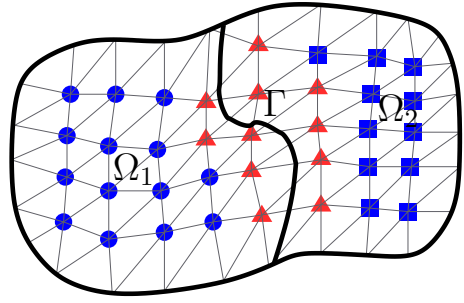}
	\put(23,26){\Large {$\Omega_1$}}
	\put(81,39){\Large {$\Omega_2$}}
	\put(56,40){\Large {$\Gamma$}}
      \end{overpic}
\caption{
DOF are associated to the interface $\Gamma$ when their basis functions have supports that are intersecting
the interface $\Gamma$.}
\label{AMLS classic: subfigure: Omega substructuring triangulation}
\end{subfigure}
\caption{
Substructuring of $\Omega$ with triangulation. Degrees of freedom (DOF) of the finite 
element space of continuous and piecewise affine functions are indicated 
by circles if they are associated to $\Omega_1$, 
by squares if associated to $\Omega_2$,
and by triangles if associated to $\Gamma$.
} 
\label{AMLS classic: figure: Omega triangulation}
\end{figure}

To solve eigenvalue problem \eqref{AMLS classic: discrete EVP}
by the AMLS method, the underlying domain $\Omega$ 
in \eqref{AMLS classic: continuous EVP classic}
is substructured into two non-overlapping 
subdomains $\Omega_1$
and $\Omega_2$ 
with
\begin{equation*}
\overline{\Omega} = \overline{\Omega}_1 \cup \overline{\Omega}_2 
\qquad
\textnormal{and}
\qquad
\Omega_1 \cap \Omega_2 = \emptyset,
\end{equation*}
which share the interface $\Gamma := \overline{\Omega}_1 \cap \overline{\Omega}_2$
(an example is given in Figure \ref{AMLS classic: subfigure: Omega substructuring}).
Since the matrices $K,M \in \R^{N \times N}$ in \eqref{AMLS classic: discrete EVP} result 
from a finite element discretisation
each row and column index is associated with
a basis function that has typically a small support.
Using the substructuring of $\Omega$ the row and column indices of $K$ and $M$ are 
reordered in such a way that it holds
\begin{equation}
K
\,
=
\bbordermatrix
{
& \Lightblue{\Omega_1} & \Lightblue{\Omega_2} & \Lightblue{\Gamma} \cr
\Lightblue{\Omega_1} & K_{11}& &K_{13} \cr
\Lightblue{\Omega_2} & & K_{22} & K_{23} \cr
\Lightblue{\Gamma} & K_{31} & K_{32} & K_{33} \cr
}
\qquad
\textnormal{and}
\qquad
M
\,
=
\bbordermatrix
{
& \Lightblue{\Omega_1} & \Lightblue{\Omega_2} & \Lightblue{\Gamma} \cr
\Lightblue{\Omega_1} & M_{11}& &M_{13} \cr
\Lightblue{\Omega_2} & & M_{22} & M_{23} \cr
\Lightblue{\Gamma} & M_{31} & M_{32} & M_{33} \cr
} 
\label{AMLS classic: matrix partitioning}
\end{equation}
where $K_{ij}, M_{ij} \in \R^{N_i \times N_j}$ and $N_1+N_2+N_3 = N$. 
The labels $\Omega_1, \Omega_2,\Gamma$ used in \eqref{AMLS classic: matrix partitioning} are indicating to which
subset the row and column indices are associated, i.e., if the supports of the corresponding basis 
functions are inside
$\Omega_i$ or if they are intersecting $\Gamma$ (cf. Figure \ref{AMLS classic: subfigure: Omega substructuring 
triangulation}).

In the next step of AMLS a block \LDL-decomposition is performed in order to 
block diagonalise the stiffness matrix $K$ by $K = L \widetilde{K} L^T$ where
\begin{equation}
L
:=
\begin{bmatrix}
\Id&       &           \\
      & \Id&           \\
K_{31} K_{11}^{-1}& K_{32} K_{22}^{-1}& \Id\\
\end{bmatrix}
\in \R^{N \times N}
\qquad
\textnormal{and}
\qquad
\widetilde{K} = \diag \bigl [{K}_{11}, {K}_{22},\widetilde{K}_{33} \bigr ],
\label{AMLS classic: block diagonalise K}
\end{equation}
and where the submatrix $\widetilde{K}_{33}$ given by
\begin{equation*}
\widetilde{K}_{33} 
= \, 
K_{33} - K_{31} K_{11}^{-1} K_{13} - K_{32} K_{22}^{-1} K_{23}
\end{equation*}
is the \emph{Schur complement} of $\diag[K_1,K_2]$ in $K$.
The mass matrix $M$ is transformed correspondingly by computing $\widetilde{M} := L^{-1} M L^{-T}$
with
\begin{equation}
\widetilde{M}
= 
\begin{bmatrix}
{M}_{11}&       & \widetilde{M}_{13}\\
      & {M}_{22}& \widetilde{M}_{23}\\
\widetilde{M}_{31}& \widetilde{M}_{32}& \widetilde{M}_{33}
\label{AMLS classic: transformation M}
\end{bmatrix},
\end{equation}
where the submatrices of $\widetilde{M}$ are given by
\begin{align*}
\widetilde{M}_{3i}
\; &= \;
M_{3i} 
-
K_{3i} K_{ii}^{-1} M_{ii}
\qquad \textnormal{and} \qquad 
\widetilde{M}_{i3} 
\; = \;
\widetilde{M}_{3i}^{T}
\qquad 
(\textnormal{for } i=1,2),
\\
\widetilde{M}_{33}
\; &= \;
M_{33} 
- 
\sum_{i=1}^{2} 
\bigl ( \,
K_{3i} K_{ii}^{-1} M_{i3}
+ 
M_{3i} K_{ii}^{-1} K_{i3}
-
K_{3i} K_{ii}^{-1} M_{ii} K_{ii}^{-1} K_{i3} 
\, \bigr ).
\end{align*}
Note that the eigenvalue problems $(\widetilde{K},\widetilde{M})$
and $(K,M)$ are equivalent, i.e.,
the eigenvalues of both problems are equal and if $\widetilde{x}$ is an eigenvector of
problem $(\widetilde{K},\widetilde{M})$ then ${x}=L^{-T}\widetilde{x}$ is an eigenvector of $(K,M)$.
The reason for the performed problem transformation from $(K,M)$ to 
$(\widetilde{K},\widetilde{M})$ is explained in the best way by describing AMLS in a continuous setting, however, 
this issue is not discussed in this paper and instead it is referred 
to, e.g., \cite{Bennighof_Lehoucq,Diss_Gerds,Quarteroni}.

In the next step of AMLS partial eigensolutions of the subproblems 
$({K}_{11},{M}_{11})$, $({K}_{22},{M}_{22})$ and $(\widetilde{K}_{33},\widetilde{M}_{33})$
are computed. To be more precise,
only those eigenpairs of the subproblems are computed
which belong to the smallest $k_i \in \N$ eigenvalues 
for given $k_i \leq N_i$ ($i=1,2,3$). In the following these partial eigensolutions are 
\begin{equation}
{K}_{ii} \, \widetilde{S}_{i} 
\; = \;
{M}_{ii} \, \widetilde{S}_{i} \, \widetilde{D}_i 
\quad \textnormal{for } i=1,2
\qquad
\textnormal{and}
\qquad
\widetilde{K}_{33} \, \widetilde{S}_{3} 
\; = \;
\widetilde{M}_{33} \, \widetilde{S}_{3} \, \widetilde{D}_3
\label{AMLS classic: partial eigensolution},
\end{equation}
where the diagonal matrix $\widetilde{D}_i \in \R^{{k_i} \times {k_i}}$ contains the $k_i$ smallest
eigenvalues, and where the matrix $\widetilde{S}_i \in \R^{N_i \times {k_i}}$ contains column-wise the associated
eigenvectors ($i=1,2,3$) which are normalised by $\widetilde{S}_i^{T} {M}_{ii} \widetilde{S}_i = \Id$ 
($i=1,2$) and $\widetilde{S}_3^{T} \widetilde{M}_{33} \widetilde{S}_3 = \Id$.

It is important to note that 
the original eigenvalue problem $(K,M)$ is not solved just by 
computing the eigenpairs of the subproblems $(K_{11},M_{11})$, $(K_{22},M_{22})$ and 
$(\widetilde{K}_{33},\widetilde{M}_{33})$. 
However, 
it could be shown (see, e.g., \cite{Bennighof_Lehoucq,Analysis_and_comparison_CMS,CMS_second_order_operators})
that the eigenvectors of the subproblems, which are associated to the smallest eigenvalues, 
are very well suited to form a subspace for the approximation of the original 
eigenvalue problem $(K,M)$.
In order to approximate eigenvalue problem $(K,M)$, respectively the sought eigenpairs of $(K,M)$, in the next step of 
AMLS the block diagonal matrix 
\begin{equation*}
\Ssum := 
\diag
\left [
\widetilde{S}_{1}, \widetilde{S}_{2}, \widetilde{S}_{3}
\right ]
\in \R^{N \times \ksum}
\qquad 
\textnormal{with } \ksum := k_1 + k_2 + k_3 \ll N
\end{equation*}
is defined, and the \emph{reduced matrices}
$
\widehat{K} := \Ssum^{T} \, \widetilde{K} \, \Ssum 
$
and
$
\widehat{M} := \Ssum^{T} \, \widetilde{M} \, \Ssum
$
are computed where it holds
\begin{equation*}
\widehat{K}
= 
\diag
\left [
\widetilde{D}_1,\widetilde{D}_2, \widetilde{D}_3\\
\right ]
\in \R^{\ksum \times \ksum}
\qquad
\textnormal{and}
\qquad
\widehat{M} = 
\begin{bmatrix}
\Id&       & \widehat{M}_{13}\\
      & \Id& \widehat{M}_{23}\\
\widehat{M}_{31}& \widehat{M}_{32}& \Id\\
\end{bmatrix}
\in \R^{\ksum \times \ksum}.
\end{equation*}
With the reduced matrices a
\emph{reduced eigenvalue problem} 
\begin{equation}
\begin{cases} 
\; \textnormal{find }  ( \, \widehat{\lambda}, \widehat{x} \, ) \in  \R \times \R^{\ksum} \setminusNull \textnormal{ 
with}
\rule[-3mm]{0mm}{3mm} 
\\
\; \widehat{K} \, \widehat{x} \, = \, \widehat{\lambda} \, \widehat{M} \, \widehat{x}
\end{cases}
\label{AMLS classic: AMLS reduced EVP}
\end{equation}
arises which possesses the eigenpairs
$$
\bigl (
\, \widehat{\lambda}_j,\widehat{x}_j
\bigr )_{j=1}^{\ksum} \in \R_{>0} \times \R^{\ksum} \setminusNull
\qquad \textnormal{and} \quad
\widehat{\lambda}_j \leq \widehat{\lambda}_{j+1}.
$$
Finally, in the last steps of AMLS the smallest $\nev$ eigenpairs of \eqref{AMLS classic: AMLS reduced EVP} are 
computed and the eigenvectors $\widehat{x}_j$ of the reduced problem are transformed by computing
\begin{equation}
\widehat{y}_j:= L^{-T}  \Ssum \; \widehat{x}_j
\qquad 
\textnormal{with }
j = 1,\ldots,\nev.
\label{AMLS classic: AMLS eigenvector approxiamtion}
\end{equation}
The vectors $\widehat{y}_j$ are Ritz-vectors of the original eigenvalue problem $(K,M)$ 
respective to the subspace spanned by the columns of the matrix $L^{-T}\Ssum$, and $\widehat{\lambda}_j$ are the
respective Ritz-values.
The Ritz-pairs
\begin{equation}
\bigl (
\, \widehat{\lambda}_j,\widehat{y}_j
\bigr )_{j=1}^{\nev} \in \R_{>0} \times \R^N \setminusNull
\qquad \textnormal{with} \quad
\widehat{\lambda}_j \leq \widehat{\lambda}_{j+1}
\label{AMLS classic: AMLS eigenpair approxiamtion before backpermutation}
\end{equation}
are approximating the sought smallest $\nev$ eigenpairs
of the eigenvalue problem $(K,M)$.

\begin{remark}
\label{AMLS classic: remark: mode selection and reduced problem}

\begin{itemize}

\item[i)]
Note that the Ritz-pairs $(\widehat{\lambda}_j,\widehat{y}_j)$, which are computed by AMLS,
are primarily used to approximate the eigensolutions of the continuous problem 
\eqref{AMLS classic: continuous EVP classic}
and not the eigenpairs of the discretised problem $(K,M)$. Correspondingly the approximation error of AMLS is
influenced by both the modal truncation performed in \eqref{AMLS classic: partial eigensolution}
and the finite element discretisation.
In contrast to the AMLS method, the approximation error of a classical approach, 
like the Lanczos method \cite{Templates_EWP} or the QR-algorithm \cite{QR},
is only influenced by the finite element discretisation
since (almost) exact eigenpairs of the discrete problem $(K,M)$ are computed.
This means as long as the error caused by the modal truncation in AMLS is of the same order as the discretisation 
error, the
eigenpair approximations computed by AMLS are of comparable quality as the eigenpair approximations computed by a 
classical approach.

\item[ii)]
How many eigenvectors have to be computed in \eqref{AMLS classic: partial eigensolution} for each subproblem is not
easy to answer. On the one hand $k_i$
should be large enough in order to keep enough spectral information from each subproblem 
so that sufficiently good eigenpair approximations are obtained from the reduced problem 
$(\widehat{K},\widehat{M})$
[e.g., AMLS computes 
exact eigenpairs of the discrete problem $(K,M)$ if $k_i=N_i$ for $i=1,2,3$]. 
On the other hand $k_i$ should be small in order to obtain a reduced problem $(\widehat{K},\widehat{M})$
of small size which can be easily solved.
Several heuristic approaches have
been discussed in literature for the eigenpair selection in \eqref{AMLS classic: partial eigensolution}:
One possible strategy is to select from each subproblem only those eigenvectors whose eigenvalues are smaller 
than a given \emph{truncation bound} $\omega {>0}$ (see, e.g., \cite{Voss,Yang_Mode_Selection}).
Another possible strategy, which is motivated by approximation properties of the finite element discretisation
of elliptic PDE eigenvalue problems, is to compute only these eigenpairs of the subproblems which 
belong, for example, to the smallest $k_i$ eigenvalues with
\begin{equation}
k_{i} = C  N_i^{1/3} \in \N
\quad
\textnormal{ or }
\quad
k_{i} = C  N_i^{1/2} \in \N
\label{AMLS classic: size of k_i mode selection}
\end{equation}
and some constant $C>0$ (see \cite{Diss_Gerds,Gerds} for details).

\item[iii)]
Since the number $k_i$ of selected eigenpairs in 
\eqref{AMLS classic: partial eigensolution} is typically quite small,
the size $\ksum$ of the reduced eigenvalue problem $(\widehat{K},\widehat{M})$
is much smaller than the size $N$ of the original problem $(K,M)$, 
and hence the reduced problem is typically much easier to solve than 
the original problem. 
If for example the mode selection described in ii)
is used with $k_i = C N_i^{1/3}$
then the size of the reduced problem can be bounded by $\mathcal{O}(N^{1/3})$
and the problem can be solved by dense linear algebra routines in ${\cal O}(N)$.

\end{itemize}

\end{remark}

\section{Introduction of the new dense AMLS method}
\label{section: AMLS dense} 
The new domain decomposition method for dense eigenvalue problems, 
which is introduced in this section and which is called in the following \emph{dense AMLS method}, 
is motivated by the classical AMLS method from the previous section. 
The initial point of dense AMLS is the finite element discretisation 
\eqref{introduction: discrete EVP} of the continuous 
integral equation eigenvalue
problem \eqref{introduction: continuous EVP classic}.
For reasons of convenience
the upper index of $\lambda\hidx, x\hidx, K\hidx$ and $M\hidx$,
which is indicating in \eqref{introduction: discrete EVP} 
the mesh width $h$ of the underlying finite element discretisation,
is left out in this particular section,
i.e., in the following we consider the eigenvalue problem 
\begin{equation}
\begin{cases} 
\; \textnormal{find } \left (\lambda, x \right ) \in  \R \times \R^N \setminusNull  \textnormal{ with}
\\
\; K \, x \, = \, \lambda \, M \, x
\end{cases}
\label{AMLS: discrete EVP}
\end{equation}
with the eigenpairs
$$
\bigl (
\lambda_j,x_j
\bigr )_{j=1}^{N} \in \R \times \R^N \setminusNull
\qquad \textnormal{where} \quad
|\lambda_j| \geq |\lambda_{j+1}|
$$
and with $N:=N_h =\dim V_h$. To avoid misunderstandings, it is explicitly noted that in this section $\lambda$ and 
$\lambda_j$ are interpreted as the eigenvalues of the discrete problem \eqref{AMLS: discrete EVP} and not as the 
eigenvalues of the continuous problem \eqref{introduction: continuous EVP classic}.
Furthermore, we are interested only in computing the by magnitude largest eigenpairs of 
\eqref{AMLS: discrete EVP},
where the number of sought eigenpairs $\nev$ is, e.g., of the size
\eqref{introduction: number of sought eigensolutions}.
\\

\begin{figure}[t]
\begin{subfigure}[t]{0.47\textwidth}
\centering
      \begin{overpic}[width=7cm]{domain_lvl_1.eps}
	\put(23,26){\Large {$\Omega_1$}}
	\put(81,39){\Large {$\Omega_2$}}
      \end{overpic}
\caption{
Substructuring of the domain $\Omega$ into the non-overlapping subdomains $\Omega_1$ and $\Omega_2$,
where $\overline{\Omega} = \overline{\Omega}_1 \cup \overline{\Omega}_2$ and 
$\Omega_1 \cap \Omega_2 = \emptyset$.}
\label{AMLS: subfigure: Omega substructuring}
\end{subfigure}
\hfill
\begin{subfigure}[t]{0.47\textwidth}
\centering
      \begin{overpic}[width=7cm]{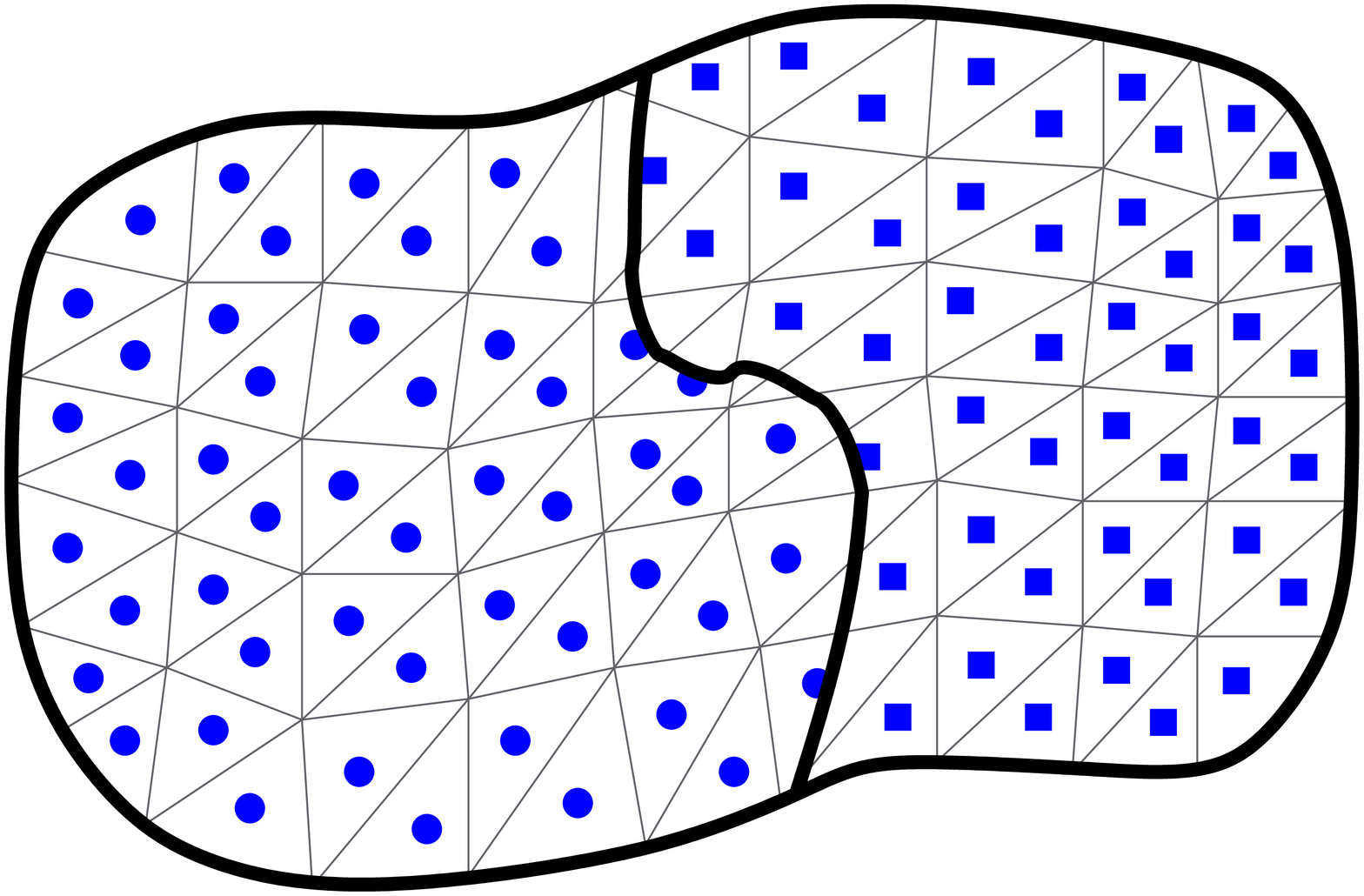}
	\put(23,26){\colorbox{white}{\Large {$\Omega_1$}}}
	\put(81,39){\colorbox{white}{\Large {$\Omega_2$}}}
      \end{overpic}
\caption{
DOF are associated to $\Omega_i$ when the nodal points of the corresponding basis functions
are elements of $\Omega_i$.}
\label{AMLS: subfigure: Omega substructuring triangulation}
\end{subfigure}
\caption{
Substructuring of $\Omega$ with triangulation. Degrees of freedom (DOF) of the finite 
element space of piecewise constant basis functions are indicated 
by circles if they are associated to $\Omega_1$ and by squares if associated to $\Omega_2$.
} 
\label{AMLS: figure: Omega triangulation}
\end{figure}

In the first step of dense AMLS the domain $\Omega$ is subdivided into two non-overlapping 
subdomains $\Omega_1$ and $\Omega_2$
(cf. Figure \ref{AMLS: subfigure: Omega substructuring})
such that
\begin{equation}
\overline{\Omega} = \overline{\Omega}_1 \cup \overline{\Omega}_2 
\qquad
\textnormal{and}
\qquad
\Omega_1 \cap \Omega_2 = \emptyset.
\label{AMLS: subdomain decomposition lvl 1}
\end{equation}
Since the matrices $K,M \in \R^{N \times N}$ in \eqref{AMLS: discrete EVP} result from a 
finite element discretisation
each row and column index is associated with
a basis function, which has typically a small support
and which has some nodal point (this might be, e.g., the midpoint of the corresponding finite element).
Using the substructuring of $\Omega$ the row and column indices of $K$ and $M$ are 
reordered in such a way that
\begin{equation}
K
\,
=
\bbordermatrix
{
& \Lightblue{\Omega_1} & \Lightblue{\Omega_2} \cr
\Lightblue{\Omega_1} & K_{11}& K_{12} \cr
\Lightblue{\Omega_2} & K_{21} & K_{22} \cr
}
\qquad
\textnormal{and}
\qquad
M
\,
=
\bbordermatrix
{
& \Lightblue{\Omega_1} & \Lightblue{\Omega_2} \cr
\Lightblue{\Omega_1} & M_{11}& M_{12} \cr
\Lightblue{\Omega_2} & M_{21} & M_{22} \cr
}
\label{AMLS: matrix partitioning}
\end{equation}
holds with $K_{ij}, M_{ij} \in \R^{N_i \times N_j}$ and $N_1+N_2 = N$. 
The labels $\Omega_1$ and $\Omega_2$ in \eqref{AMLS: matrix partitioning} are indicating to which
subset the row and column indices are associated.
In contrast to the matrix partitioning performed in 
classical AMLS (which depends on the supports of the basis functions),
the matrix partitioning in dense AMLS depends only on the nodal points of the corresponding basis functions.
A row or column index is associated to $\Omega_1$ if the nodal point of the corresponding basis function is inside
$\Omega_1$, otherwise the index is associated to $\Omega_2$
(cf. Figure \ref{AMLS: subfigure: Omega substructuring triangulation}).

In the next step of dense AMLS a block \LDL-decomposition is performed in order to 
block diagonalise the matrix $K$ by $K = L \widetilde{K} L^T$ where
\begin{equation}
L
:=
\begin{bmatrix}
\Id&                \\
K_{21} K_{11}^{-1} & \Id\\
\end{bmatrix}
\in \R^{N \times N}
\qquad
\textnormal{and}
\qquad
\widetilde{K} = \diag \bigl [{K}_{11}, \widetilde{K}_{22} \bigr ].
\label{AMLS: block diagonalise K}
\end{equation}
The submatrix $\widetilde{K}_{33}$ given by
\begin{equation*}
\widetilde{K}_{22} 
= \,
K_{22} - K_{21} K_{11}^{-1} K_{12}
\end{equation*}
is the \emph{Schur complement} of $K_{11}$ in $K$.
The matrix $M$ is transformed correspondingly by computing $\widetilde{M} := L^{-1} M L^{-T}$
with
\begin{equation}
\widetilde{M}
= 
\begin{bmatrix}
{M}_{11}&    \widetilde{M}_{12}\\
\widetilde{M}_{21}& \widetilde{M}_{22}
\label{AMLS: transformation M}
\end{bmatrix}
\end{equation}
where the submatrices of $\widetilde{M}$ are given by
\begin{align*}
\widetilde{M}_{21}
\; &= \;
M_{21} 
-
K_{21} K_{11}^{-1} M_{11}
\qquad \textnormal{and} \qquad 
\widetilde{M}_{12} 
\; = \;
\widetilde{M}_{21}^{T},
\\
\widetilde{M}_{22}
\; &= \;
M_{22} 
-
K_{21} K_{11}^{-1} M_{12}
-
M_{21} K_{11}^{-1} K_{12}
-
K_{21} K_{11}^{-1} M_{11} K_{11}^{-1} K_{12}.
\end{align*}
The eigenvalue problems $(\widetilde{K},\widetilde{M})$ and $(K,M)$ are equivalent, i.e.,
the eigenvalues of both problems are equal and if $\widetilde{x}$ is an eigenvector of
$(\widetilde{K},\widetilde{M})$ then ${x}=L^{-T}\widetilde{x}$ is an eigenvector of $(K,M)$.

In the next step of dense AMLS partial eigensolutions of the subproblems 
$({K}_{11},{M}_{11})$ and $(\widetilde{K}_{22},\widetilde{M}_{22})$
are computed.
However,
in contrast to the classical AMLS method, only those eigenpairs of the subproblems
are computed which belong to the by magnitude
largest 
$k_i \in \N$ eigenvalues for given $k_i \leq N_i$ ($i=1,2$). In the following these partial eigensolutions are 
\begin{equation}
{K}_{11} \, \widetilde{S}_{1} 
\; = \;
{M}_{11} \, \widetilde{S}_{1} \, \widetilde{D}_1
\qquad
\textnormal{and}
\qquad
\widetilde{K}_{22} \, \widetilde{S}_{2} 
\; = \;
\widetilde{M}_{22} \, \widetilde{S}_{2} \, \widetilde{D}_2,
\label{AMLS: partial eigensolution}
\end{equation}
where the diagonal matrix $\widetilde{D}_i \in \R^{{k_i} \times {k_i}}$ contains the $k_i$ selected
eigenvalues, and the matrix $\widetilde{S}_i \in \R^{N_i \times {k_i}}$ contains column-wise the associated
eigenvectors ($i=1,2$) which are normalised by 
$\widetilde{S}_1^{T} {M}_{11} \widetilde{S}_1 = \Id$ 
and
$\widetilde{S}_2^{T} \widetilde{M}_{22} \widetilde{S}_2 = \Id$.
\begin{remark}[Mode Selection]
\label{AMLS: remark: mode selection}
Motivated by the mode selection strategy
described in Remark 
\ref{AMLS classic: remark: mode selection and reduced problem}\:ii), it is proposed to 
compute in the dense AMLS method (which has been introduced in the framework of 
integral equation eigenvalue problems) only these eigenpairs of 
\eqref{AMLS: partial eigensolution} which belong to the by magnitude largest $k_i$ eigenvalues and where 
$k_i$ is chosen, for example, as in \eqref{AMLS classic: size of k_i mode selection}.
\end{remark}

In the next step of dense AMLS we define the block diagonal matrix 
\begin{equation*}
\Ssum := 
\diag
\left [
\widetilde{S}_{1}, \widetilde{S}_{2}
\right ]
\in \R^{N \times \ksum}
\qquad 
\textnormal{with } \ksum := k_1 + k_2 \ll N.
\end{equation*} 
The $\ksum$-dimensional subspace which is spanned by the columns of the matrix $L^{-T}\Ssum$ is potentially well suited 
to approximate the sought eigenpairs of  $(K,M)$. 
To approximate these eigenpairs, in the next step of dense AMLS the \emph{reduced matrices}
$
\widehat{K} := \Ssum^{T} \, \widetilde{K} \, \Ssum 
$
and
$
\widehat{M} := \Ssum^{T} \, \widetilde{M} \, \Ssum
$
are computed where it holds
\begin{equation*}
\widehat{K}
= 
\diag
\left [
\widetilde{D}_1,\widetilde{D}_2\\
\right ]
\in \R^{\ksum \times \ksum}
\qquad
\textnormal{and}
\qquad
\widehat{M} = 
\begin{bmatrix}
\Id       & \widehat{M}_{12}\\
\widehat{M}_{21} & \Id\\
\end{bmatrix}
\in \R^{\ksum \times \ksum},
\end{equation*}
and a \emph{reduced eigenvalue problem} 
\begin{equation}
\begin{cases} 
\; \textnormal{find }  ( \, \widehat{\lambda}, \widehat{x} \, ) \in  \R \times \R^{\ksum} \setminusNull \textnormal{ 
with}
\rule[-3mm]{0mm}{3mm} 
\\
\; \widehat{K} \, \widehat{x} \, = \, \widehat{\lambda} \, \widehat{M} \, \widehat{x}
\end{cases}
\label{AMLS: AMLS reduced EVP}
\end{equation}
is obtained which possesses the eigenpairs
$$
\bigl (
\, \widehat{\lambda}_j,\widehat{x}_j
\bigr )_{j=1}^{\ksum} \in \R \times \R^{\ksum} \setminusNull
\qquad \textnormal{and} \quad
|\widehat{\lambda}_j| \geq |\widehat{\lambda}_{j+1}|.
$$

At the end of dense AMLS the by magnitude largest $\nev$ eigenpairs of \eqref{AMLS: AMLS reduced EVP} are computed
and the eigenvectors $\widehat{x}_j$ of the reduced problem are transformed via
\begin{equation}
\widehat{y}_j:= L^{-T}  \Ssum \; \widehat{x}_j
\qquad 
\textnormal{with }
j = 1,\ldots,\nev.
\label{AMLS: AMLS eigenvector approxiamtion}
\end{equation}
The vectors $\widehat{y}_j$ are the Ritz-vectors of the original eigenvalue problem $(K,M)$ 
respective to the subspace spanned by the columns of the matrix $L^{-T}\Ssum$, and $\widehat{\lambda}_j$ are the
respective Ritz-values. The Ritz-pairs 
\begin{equation}
\bigl (
\, \widehat{\lambda}_j,\widehat{y}_j
\bigr )_{j=1}^{\nev} \in \R \times \R^N \setminusNull
\qquad \textnormal{with} \quad
|\widehat{\lambda}_j| \geq |\widehat{\lambda}_{j+1}|
\label{AMLS: AMLS eigenpair approxiamtion before backpermutation}
\end{equation}
are approximating the sought by magnitude largest $\nev$ eigenpairs
of the problem $(K,M)$.

To summarise dense AMLS an overview of all necessary operations is given in
Table \ref{AMLS: table: overview AMLS} where the different tasks of the method are denoted by
$(\task 1)$--$(\task 8)$.

\begin{table}[t]
\centering
\scalebox{0.90}{
\renewcommand{\arraystretch}{2.0}
    \begin{tabular}{| p{0.6cm} p{5.7cm} | p{9.3cm} |}
    \hline
    \multicolumn{2}{|c|}{\rule[-4mm]{0mm}{10mm} task} 
    & \multicolumn{1}{c|}{matrix operations of dense AMLS}
    \\ \hline \hline 
    %
    %
    ($\task 1$) & partition matrices $K$ and $M$ &
    apply geometric bisection reordering as done in \eqref{AMLS: matrix partitioning}
    \\ \hline
    %
    %
    ($\task 2$) & block diagonalise the matrix $K$ & 
    $K = L \, \widetilde{K} \, L^T$
    \\ \hline 
    %
    %
    ($\task 3$) & transform $M$ & 
    $\widetilde{M} = L^{-1} \,  M  \, L^{-T}$
    \\ \hline 
    %
    %
    ($\task 4$) & compute partial eigensolutions \newline of the subproblems
    & 
    ${K}_{11} \, \widetilde{S}_{1} = {M}_{11} \, \widetilde{S}_{1} \, \widetilde{D}_1$
    \quad\; with $\widetilde{S}_1 \in \R^{N_1 \times k_1}$, $\widetilde{D}_1 \in \R^{k_1 \times 
k_1}$\rule[-2mm]{0mm}{0mm}
    \newline
    $\widetilde{K}_{22} \, \widetilde{S}_{2} = \widetilde{M}_{22} \, \widetilde{S}_{2} \, \widetilde{D}_2$
    \quad\; with $\widetilde{S}_2 \in \R^{N_2 \times k_2}$, $\widetilde{D}_2 \in \R^{k_2 \times k_2}$
    \\ \hline 
    %
    %
    ($\task 5$) & define subspace &
    $ \Ssum  := \, \diag
    \left [ \widetilde{S}_1, \widetilde{S}_2 \right ] \in \R^{N \times \ksum}
    \qquad \textnormal{with } \ksum := k_1+k_2$
    \\ \hline 
    %
    %
    ($\task 6$) & compute the reduced matrices & 
    $\widehat{K} := {\Ssum}^T \, \widetilde{K} \, {\Ssum} \in \R^{\ksum \times \ksum}$
    \quad and \quad
    $\widehat{M} := {\Ssum}^T \, \widetilde{M} \, {\Ssum} \in \R^{\ksum \times \ksum}$
    \\ \hline 
    %
    %
    ($\task 7$) & solve reduced eigenvalue problem &
    $\widehat{K} \, \widehat{x}_j \, = \, \widehat{\lambda}_j \, \widehat{M} \, \widehat{x}_j$
    \qquad for $j=1, \ldots, \nev$
    \\ \hline 
    %
    %
    ($\task 8$) &transform eigenvectors &
    $\widehat{y}_j:= L^{-T}  {\Ssum} \; \widehat{x}_j$ 
    \qquad for $j=1, \ldots, \nev$
    \rule[-4mm]{0mm}{0mm} 
    \\ \hline
    %
    %
    \end{tabular}
  }
  \caption{Overview of the new dense AMLS method to compute eigenpair approximations
  $( \, \widehat{\lambda}_j,\widehat{y}_j)$ for the problem $(K,M)$.
  }
  \label{AMLS: table: overview AMLS}
\end{table}

\begin{remark}[Rayleigh-Ritz Projection]
\label{AMLS: remark: Rayleigh-Ritz Projection}
\begin{itemize}

\item[i)]
The dense AMLS method can be summarised also as follows:
At first the subspace spanned by the columns of the matrix $Q:= L^{-T}  \Ssum$
is computed [task $(\task1)$ -- $(\task4)$]. Depending on the mode selection strategy in task $(\task4)$
this subspace is well suited to approximate the sought eigensolutions of $(K,M)$.
Thereafter, the Rayleigh-Ritz projection of problem $(K,M)$ onto the subspace $Q$ is computed,
i.e., the reduced matrices $\widehat{K} = Q^T K Q$ and $\widehat{M} = Q^T K Q$ are computed [task $(\task6)$], the 
reduced eigenvalue problem $(\widehat{K},\widehat{M})$ is solved [task $(\task7)$], where the pairs
$(\,\widehat{\lambda}_j,Q \widehat{x}_j)$ [task $(\task8)$]
finally provide approximations of the sought eigenpairs of $(K,M)$.

\item[ii)]
To solve eigenvalue problem \eqref{introduction: continuous EVP classic}
by dense AMLS the role of $K$ and $M$ in \eqref{AMLS: discrete EVP} can be interchanged
so that problem $Mx=1/\lambda Kx$ is considered instead.
In this setting the eigenpairs of $(M,K)$ associated to the by magnitude smallest eigenvalues are sought, and in 
task $(\task 4)$ the by magnitude smallest $k_i$ eigenpairs have to be computed.

\end{itemize}

\end{remark}

\section{Numerical Results}
\label{section: Numerical Results}

The dense AMLS method has been implemented in C++ using the LAPACK/BLAS \cite{Lapack_User_Guide,BLAS} library. In 
the following we analyse numerically the dense 
AMLS method for the integral equation eigenvalue problem
\begin{equation}
\int_{\Omega}
\log |x-y| u(y) \dy
=
\lambda u(x)
\qquad \textnormal{with } x \in \Omega:=(0,1).
\label{results: eigenvalue problem example}
\end{equation}

To solve problem \eqref{results: eigenvalue problem example} by the dense AMLS method
or by a classical approach it is discretised by the finite element method
using piecewise constant basis functions
\begin{equation}
\varphi_i\Nidx:=
\begin{cases}
1\quad \textnormal{if } \tfrac{i-1}{N} \leq x \leq \tfrac{i}{N} \\
0 \quad \textnormal{otherwise}
\end{cases}
\quad \textnormal{for } 
i=1,\ldots,N
\end{equation}
where the interval $[0,1]$ is decomposed into $N$ equispaced subintervals with mesh width $h:=1/N$.
As described in \eqref{introduction: discrete EVP} the discretisation results in the algebraic
eigenvalue problem $(K\hidx,M\hidx)$ whose discrete eigenpairs $(\lambda\hidx,x\hidx)$ are approximating
the sought eigensolutions $(\lambda,u)$ of the continuous eigenvalue problem 
\eqref{results: eigenvalue problem example}.

In the benchmarks presented in the following the eigenvalues $\lambda\hidx$ of the algebraic eigenvalue 
$(K\hidx,M\hidx)$ have been all positive, and thus the discrete eigenvalues 
$\displaystyle \lambda\hidx_j:=\lambdaplus\hidx_j$ with $\lambda\hidx_1 \geq \lambda\hidx_2 \geq \ldots \geq 0$
are approximating the positive eigenvalues
$\displaystyle \lambda_j:=\lambdaplus_j$ 
of \eqref{results: eigenvalue problem example}
with $\lambda_1 \geq \lambda_2 \geq \ldots \geq 0$.

To evaluate the approximation accuracy of the eigenvalues
$\widehat{\lambda}\hidx_j$ computed by dense AMLS
(the upper index of $\widehat{\lambda}\hidx_j$ indicates 
the mesh width $h$ of the underlying finite element discretisation),
we compare the approximation quality of dense AMLS 
[cf. Remark \ref{AMLS classic: remark: mode selection and reduced problem}\:i)]
with the approximation quality of a
classical approach by examining the relative errors
\begin{equation}
\deltaHAMLSBEM
\hspace*{+1mm}
:=
\hspace*{-2mm}
\underbrace{
\dfrac{| \, \lambda_j \, - \, \widehat{\lambda}\hidx_j\, |}{
|\lambda_j| \rule[-3mm]{0mm}{7mm} 
}}_{
\begin{tabular}{c}
\begin{footnotesize} relative error of \end{footnotesize}
\\
\begin{footnotesize} dense AMLS \end{footnotesize}
\end{tabular}
}
\hspace*{3mm} \textnormal{and} \hspace*{8mm}
\deltaRef
\hspace*{+1mm}
:=
\hspace*{-2mm}
\underbrace{
\dfrac{| \, \lambda_j \, - \,  \lambda_j\hidx \, |}{
|\lambda_j| \rule[-3mm]{0mm}{7mm} 
}}_{
\begin{tabular}{c}
\begin{footnotesize} relative error \end{footnotesize}
\\
\begin{footnotesize} of discretisation \end{footnotesize}
\end{tabular}
}
\hspace*{-3mm}
.
\label{Results: relative errors}
\end{equation}
However, since for 
problem \eqref{results: eigenvalue problem example}
the exact eigenvalues $\lambda_j$ are not known, the eigenvalues
$\lambda_j$ in 
\eqref{Results: relative errors}
are approximated by the discrete eigenvalues $\lambda_j^{(h_0)}$ which are associated to the 
very fine mesh width $h_0$:=2e-4 ($h_0$ corresponds to a discretisation with 5,000 DOF).
Using these approximated relative errors the approximation quality of dense AMLS is benchmarked for 
problem 
\eqref{results: eigenvalue problem example}
by investigating 
\begin{equation}
\gamma_{\nev}^{(h)}
:=
\max 
\left \{
\deltaHAMLSBEM \, /\; \deltaRef
\; : \;
j=1,\ldots,\nev
\right \}
\label{results: error ratios}
\end{equation}
which is the maximal ratio between the relative errors $\deltaHAMLSBEM$ of dense AMLS
and the relative discretisation errors $\deltaRef$.
If it holds
\begin{equation}
\gamma_{\nev}^{(h)} < 3
\label{Results: postulation}
\end{equation}
it can be said that the approximation error of dense AMLS is of the same order as the discretisation error,
and thus in this case the approximation quality of dense AMLS competes 
with the approximation quality of a classical approach
[cf. Remark \ref{AMLS classic: remark: mode selection and reduced problem}\:i)].

\begin{figure*}[p]
\centering
\begin{minipage}[c]{0.45\textwidth}
      \begin{overpic}[width=8.2cm]{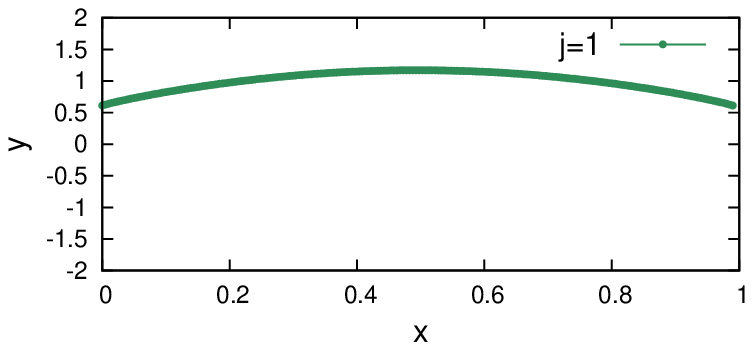}
      \end{overpic}
\end{minipage}
\hfill
\begin{minipage}[c]{0.45\textwidth}
      \begin{overpic}[width=8.2cm]{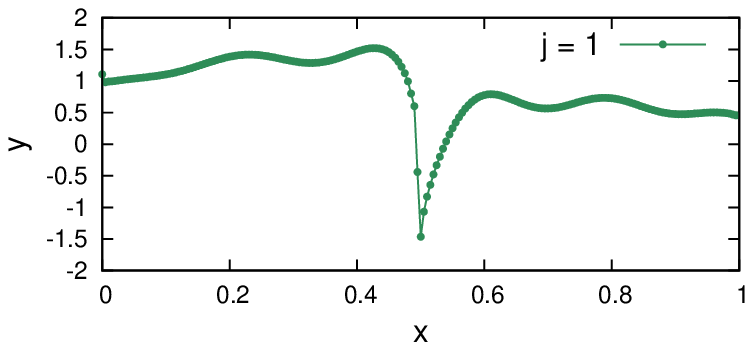}
      \end{overpic}
\end{minipage}
\\
\begin{minipage}[c]{0.45\textwidth}
      \begin{overpic}[width=8.2cm]{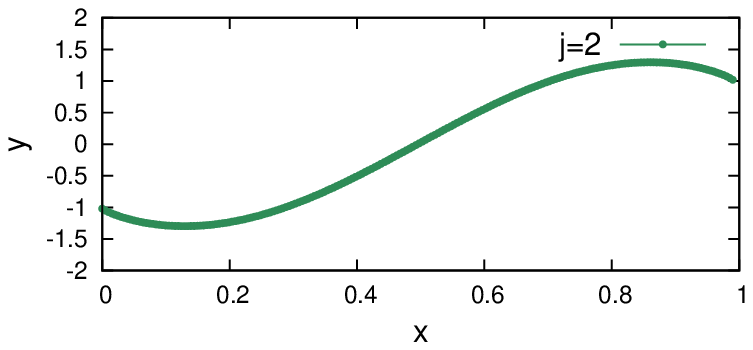}
      \end{overpic}
\end{minipage}
\hfill
\begin{minipage}[c]{0.45\textwidth}
      \begin{overpic}[width=8.2cm]{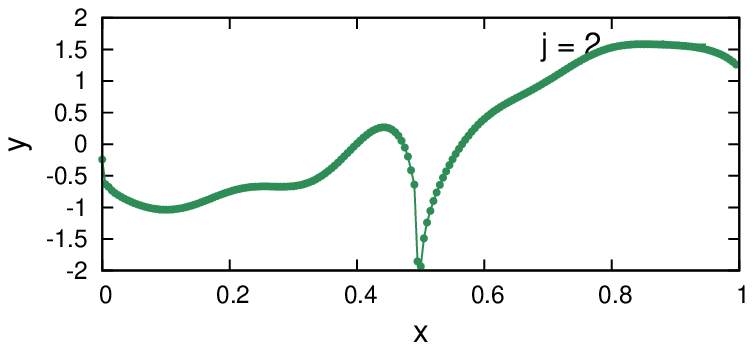}
      \end{overpic}
\end{minipage}
\\
\begin{minipage}[c]{0.45\textwidth}
      \begin{overpic}[width=8.2cm]{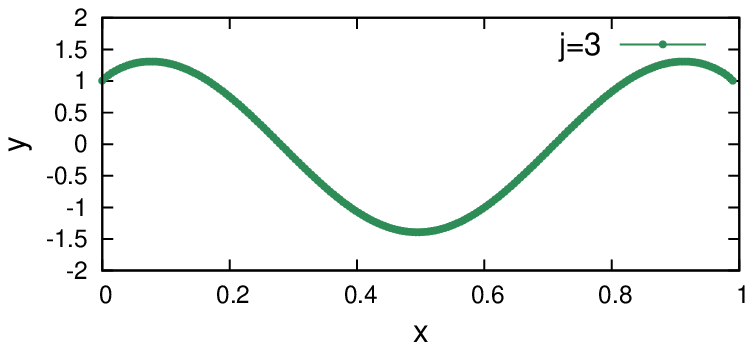}
      \end{overpic}
\end{minipage}
\hfill
\begin{minipage}[c]{0.45\textwidth}
      \begin{overpic}[width=8.2cm]{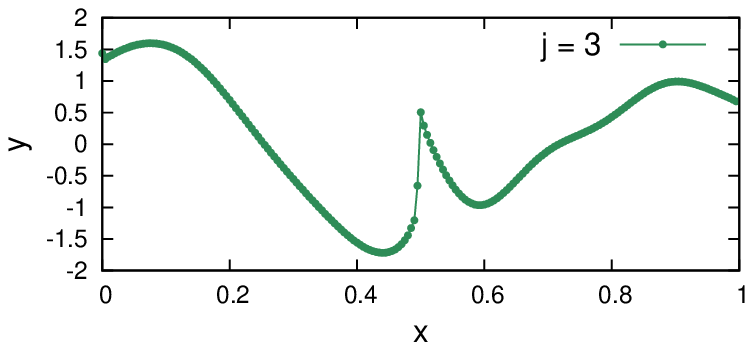}
      \end{overpic}
\end{minipage}
\\
\begin{minipage}[c]{0.45\textwidth}
      \begin{overpic}[width=8.2cm]{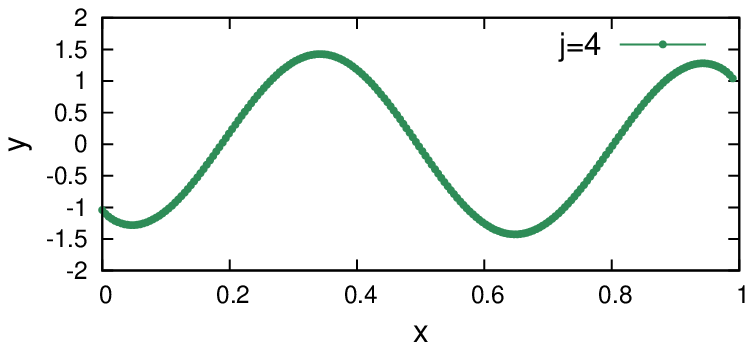}
      \end{overpic}
\end{minipage}
\hfill
\begin{minipage}[c]{0.45\textwidth}
      \begin{overpic}[width=8.2cm]{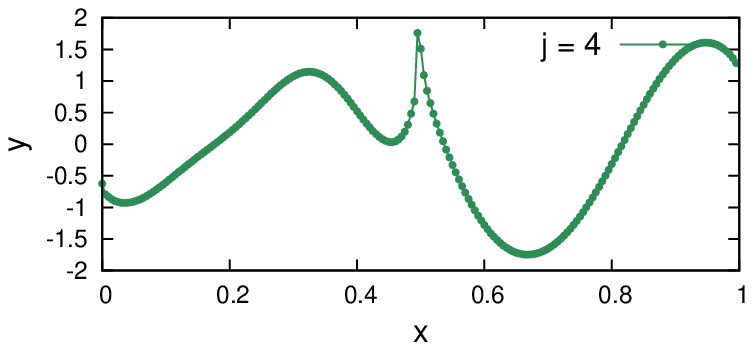}
      \end{overpic}
\end{minipage}
\\
\begin{minipage}[c]{0.45\textwidth}
      \begin{overpic}[width=8.2cm]{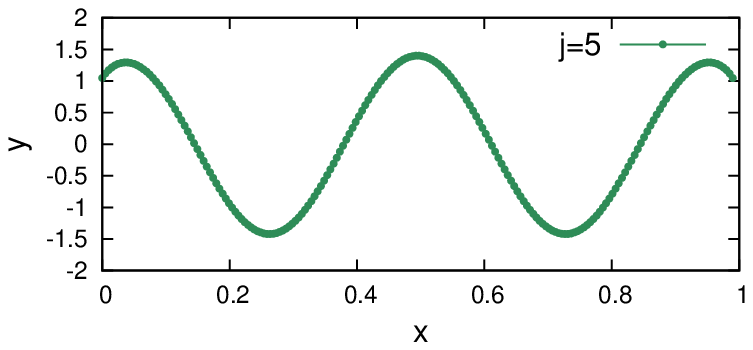}
      \end{overpic}
\end{minipage}
\hfill
\begin{minipage}[c]{0.45\textwidth}
      \begin{overpic}[width=8.2cm]{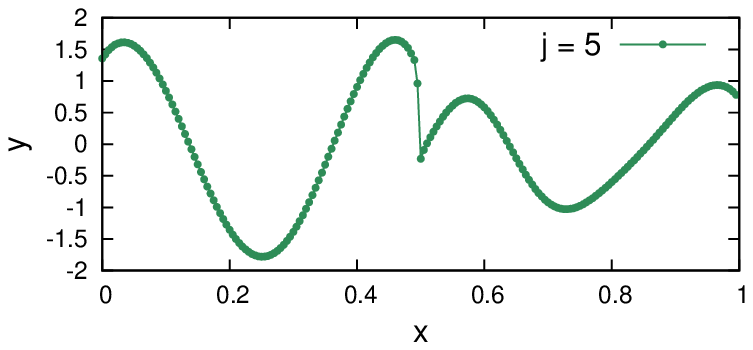}
      \end{overpic}
\end{minipage}
\caption{
Approximations for the eigenfunctions
$u_j$ of problem \eqref{results: eigenvalue problem example}
for $j=1,\ldots,5$:
On the left-hand side are the approximations provided by the discrete 
eigenvectors $x_j\hidx$ of $(K\hidx,M\hidx)$
where the underlying mesh width $h$ leads to a 
discretisation with 200 DOF, and 
on the right-hand side are the approximations that are provided by the Ritz-vectors 
$\widehat{y}_j\hidx$ of the dense AMLS method.
The function values of the eigenfunction approximations at the nodal points,
which are provided by $x_j\hidx$ and $\widehat{y}_j\hidx$,
have been connected by lines for a better visibility.}
\label{introduction: eigenfunctions classical dense AMLS}
\end{figure*}

The numerical results which are presented in the following have been obtained using the following setting:
\begin{itemize}
\item 
Eigenvalue problem \eqref{results: eigenvalue problem example} has been discretised with 
the mesh width $h$:=5e-3 leading to a discrete problem $(K\hidx,M\hidx)$ with $N_h=200$ DOF.

\item
The performed substructuring in task (\task1) of dense AMLS resulted in two subproblems 
with $100$ DOF each.
In task (\task4) in each subproblem the eigensolutions associated to the largest $k_i=5$ eigenvalues
have been computed leading to a reduced problem of size $\ksum=10$.
Hence, the dense AMLS method computed the eigenvalue approximations $\widehat{\lambda}\hidx_j$ with $j=1,\ldots,10$.
\end{itemize}
The resulting approximation accuracy of dense AMLS using the setting described above can be investigated 
in Table \ref{AMLS: table: dense AMLS}
which shows the relative errors from \eqref{Results: relative errors} and the corresponding error ratios.
The approximation accuracy of dense AMLS in this configuration is very bad. The comparison
---
see Figure \ref{introduction: eigenfunctions classical dense AMLS}
where the approximated eigenfunctions associated to the 5 largest eigenvalues are displayed
---
between the eigenfunction approximations provided by $x_j\hidx$ of the discrete problem $(K\hidx,M\hidx)$
and the eigenfunction approximations which are
provided by dense AMLS via
$\widehat{y}_j\hidx$ in  \eqref{AMLS: AMLS eigenpair approxiamtion before backpermutation}
gives a possible explanation for this unsatisfying results.
It seems that the eigenfunction approximations which are computed by dense AMLS 
are strongly varying from the exact ones especially at the region close to the interface of the subdomains 
$\Omega_1=(0,\tfrac{1}{2})$ and $\Omega_2=(\tfrac{1}{2},1)$. 
Further benchmarks showed that the approximation accuracy of dense AMLS could not be significantly 
improved just by increasing the number $k_i$ of selected eigenpairs in task (\task4).
Rather it seems that essential spectral information of the subproblems 
is missing in order to well approximate the 
sought eigenfunctions of the global problem $(K\hidx,M\hidx)$.

\begin{table*}[t]
\centering
\scalebox{0.97}{
\renewcommand{\arraystretch}{1.0}
\begin{tabular}{|r|| c|c|c|}
\hline
$j$ &  $\deltaHAMLSBEM$ & $\deltaRef$ & $\deltaHAMLSBEM/ \deltaRef$\rule[-3mm]{0mm}{10mm}
\\ \hline
\rule{0mm}{6mm} 1 & 1.93e-1 & 3.67e-6 & 5.25e+4 \\
2& 9.41e-2 & 2.74e-5 & 3.43e+3\\
3&7.72e-2 &9.70e-5  &7.95e+2\\
4&5.74e-2&2.02e-4&2.84e+2\\
5&5.10e-2&3.52e-4&1.45e+2\\
6&4.22e-2&5.38e-4&7.84e+1\\
7&4.02e-2&7.68e-4&5.23e+1\\
8&3.77e-2&1.03e-3&3.65e+1\\
9&3.67e-2&1.34e-3&2.73e+1\\
\rule[-3mm]{0mm}{0mm}10&4.64e-2&1.68e-3&2.75e+1\\
\hline
\end{tabular}
}
\caption{
Relative errors of dense AMLS 
$\deltaHAMLSBEM$
versus
relative discretisation errors
$\deltaRef$
for the eigenvalue approximation.
}
\label{AMLS: table: dense AMLS}
\end{table*}

\begin{table}[t]
\centering
\scalebox{0.90}{
\renewcommand{\arraystretch}{2.0}
    \begin{tabular}{| p{0.8cm} p{4cm} | p{4.7cm} | p{4.7cm} | }
    \hline
    \multicolumn{2}{|c|}{\rule[-4mm]{0mm}{10mm} task} 
    & \multicolumn{2}{c|}{matrix operations of combined dense AMLS}
    \\ \hline \hline 
    %
    %
    ($\ctask 1$) & partition matrices \newline $K$ and $M$ &
    apply geometric bisection \newline reordering as done in \eqref{AMLS: matrix partitioning} &
    apply geometric bisection \newline reordering as done in \eqref{AMLS: matrix partitioning changed}
    \\ \hline
    %
    %
    ($\ctask 2$) & block diagonalise $K$ & 
    $K = L \, \widetilde{K} \, L^T$ &
    $K = L \, \widetilde{K} \, L^T$
    \\ \hline 
    %
    %
    ($\ctask 3$) & transform $M$ & 
    $\widetilde{M} = L^{-1} \,  M  \, L^{-T}$ &
    $\widetilde{M} = L^{-1} \,  M  \, L^{-T}$
    \\ \hline 
    %
    %
    ($\ctask 4$) & compute partial \newline eigensolutions 
    & 
    ${K}_{11} \, \widetilde{S}_{1} = {M}_{11} \, \widetilde{S}_{1} \, \widetilde{D}_1$ \rule[-2mm]{0mm}{0mm}
    \newline
    $\widetilde{K}_{22} \, \widetilde{S}_{2} = \widetilde{M}_{22} \, \widetilde{S}_{2} \, \widetilde{D}_2$ &
    ${K}_{11} \, \widetilde{S}_{1} = {M}_{11} \, \widetilde{S}_{1} \, \widetilde{D}_1$ \rule[-2mm]{0mm}{0mm}
    \newline
    $\widetilde{K}_{22} \, \widetilde{S}_{2} = \widetilde{M}_{22} \, \widetilde{S}_{2} \, \widetilde{D}_2$ 
    \\ \hline 
    %
    %
    ($\ctask 5$) & a) compute subspaces &
    $ \Ssum  := \, \diag
    \left [ \widetilde{S}_1, \widetilde{S}_2 \right ] 
    $ \rule[-4mm]{0mm}{0mm}
    \newline 
    $Q_A:=L^{-T}Z  \in \R^{N \times \ksum_A}$
    &
    $ \Ssum  := \, \diag
    \left [ \widetilde{S}_1, \widetilde{S}_2 \right ] 
    $ \rule[-4mm]{0mm}{0mm}
    \newline 
    $Q_B:=L^{-T}Z \in \R^{N \times \ksum_B}$
    \\ \cdashline{2-4}
    %
    %
    & b) combine subspaces &
    \multicolumn{2}{|l|}{
    $ Q := \,    \left [ Q_A,Q_B \right ] \in \R^{N \times \ksum}$ \quad with $\ksum \leq \ksum_A + \ksum_B$
    }
    \\
    & &
    \multicolumn{2}{|l|}{
    (note: linear dependent columns in $Q$ have to be removed)
    }
    \\ \hline 
    %
    %
    ($\ctask 6$) & compute the \newline reduced matrices & 
    \multicolumn{2}{|l|}{
    $\widehat{K} := Q^T \, K \, Q  \in \R^{\ksum \times \ksum}$ \quad and \quad 
    $\widehat{M} := Q^T \, M \, Q \in \R^{\ksum \times \ksum}$
    }
    \\ \hline 
    %
    %
    ($\ctask 7$) & solve the reduced \newline eigenvalue problem &
    \multicolumn{2}{|l|}{
    $\widehat{K} \, \widehat{x}_j \, = \, \widehat{\lambda}_j \, \widehat{M} \, \widehat{x}_j$
    \qquad for $j=1, \ldots, \nev$
    }
    \\ \hline 
    %
    %
    ($\ctask 8$) &transform eigenvectors &
     \multicolumn{2}{|l|}{
    $\widehat{y}_j:= Q \; \widehat{x}_j$ 
    \qquad for $j=1, \ldots, \nev$
    \rule[-4mm]{0mm}{0mm} \newline
    }
    \\ \hline
    %
    %
    \end{tabular}
  }
  \caption{Overview of the combined dense AMLS method to compute eigenpair approximations
  $( \, \widehat{\lambda}_j,\widehat{y}_j)$ for the problem $(K,M)$.
  }
  \label{AMLS: table: overview combined AMLS}
\end{table}

\begin{table*}[t]
\centering
\scalebox{0.97}{
\renewcommand{\arraystretch}{1.0}
\begin{tabular}{|r|| c|c|c|}
\hline
$j$ &  $\deltaHAMLSBEM$ & $\deltaRef$ & $\deltaHAMLSBEM/ \deltaRef$\rule[-3mm]{0mm}{10mm}
\\ \hline 
\rule[0mm]{0mm}{6mm}
	1	 &9.85e-6	 &3.67e-6	 &2.68e+0\\
	 2	 &2.89e-5	 &2.74e-5	&1.06e+0\\
	 3	 &1.08e-4	 &9.70e-5	&1.11e+0\\
	 4	 &2.12e-4	 &2.02e-4	&1.05e+0\\
	 5	 &3.79e-4	 &3.52e-4	&1.08e+0\\
	 6	 &5.44e-4	 &5.38e-4	&1.01e+0\\
	 7	 &7.94e-4	 &7.68e-4	&1.03e+0\\
	 8	 &1.05e-3	 &1.03e-3	&1.01e+0\\
	 9	 &1.38e-3	 &1.34e-3	&1.03e+0\\
\rule[-2mm]{0mm}{0mm}	 10	 &1.69e-3	 &1.68e-3	&1.00e+0\\
	 \hline
\rule[0mm]{0mm}{4mm}	 11	 &2.23e-3	 &2.07e-3	&1.08e+0\\
	 12	 &5.05e-3	 &2.49e-3	&2.03e+0\\
	 13	 &3.16e-2	 &2.95e-3	&1.07e+1\\
	 14	 &9.32e-2	 &3.45e-3	&2.70e+1\\
	 15	 &2.09e-1	 &3.99e-3	&5.24e+1\\
	 16	 &4.91e-1	 &4.56e-3	&1.08e+2\\
	 17	 &4.80e-1	 &5.17e-3	&9.29e+1\\
	 18	 &5.37e-1	 &5.82e-3	&9.22e+1\\
	 19	 &6.09e-1	 &6.50e-3	&9.36e+1\\
\rule[-3mm]{0mm}{0mm}
	 20	 &8.84e-1	 &7.22e-3	&1.22e+2\\
\hline
\end{tabular}
}
\caption{
Relative errors of combined dense AMLS 
$\deltaHAMLSBEM$
versus
relative discretisation errors
$\deltaRef$
for the eigenvalue approximation.
}
\label{AMLS: table: combined dense AMLS}
\end{table*}

Furthermore, it is noted that 
in general the eigenpair approximations, which are computed by dense AMLS, differ from each other if 
in dense AMLS instead of 
\eqref{AMLS: matrix partitioning}
the changed matrix partitioning 
\begin{equation}
K
\,
=
\bbordermatrix
{
& \Lightblue{\Omega_2} & \Lightblue{\Omega_1} \cr
\Lightblue{\Omega_2} & K_{22}& K_{12} \cr
\Lightblue{\Omega_1} & K_{21} & K_{11} \cr
}
\qquad
\textnormal{and}
\qquad
M
\,
=
\bbordermatrix
{
& \Lightblue{\Omega_2} & \Lightblue{\Omega_1} \cr
\Lightblue{\Omega_2} & M_{22}& M_{12} \cr
\Lightblue{\Omega_1} & M_{21} & M_{11} \cr
}
\label{AMLS: matrix partitioning changed}
\end{equation}
is applied.
The reason for this, in simple terms, is that when matrix partitioning 
\eqref{AMLS: matrix partitioning} is used then spectral information is transferred 
from $\Omega_1$
to the subproblem 
associated to subdomain $\Omega_2$, and when matrix partitioning 
\eqref{AMLS: matrix partitioning changed} is used instead
then spectral information is transferred 
from $\Omega_2$
to the subproblem 
associated to $\Omega_1$.

As already noted in Remark \ref{AMLS: remark: Rayleigh-Ritz Projection}\:i)
the dense AMLS method can be interpreted as a Rayleigh-Ritz projection using the approximation subspace
which is spanned by the columns of the matrix $Q:=L^{-T}Z$.
This interpretation of dense AMLS motivates to apply a Rayleigh-Ritz projection for the eigenvalue problem $(K,M)$ where
the associated approximation subspace is spanned by the columns of the matrix $Q_A:=L^{-T}Z$
which is obtained by applying tasks (\task1)--(\task4) 
of dense AMLS using matrix partitioning \eqref{AMLS: matrix partitioning}
and by the columns of the matrix $Q_B:=L^{-T}Z$
which is obtained by applying tasks (\task1)--(\task4) of 
dense AMLS using matrix partitioning \eqref{AMLS: matrix partitioning changed}.
In the following this approach is referred to
as \emph{combined dense AMLS}. An overview of all necessary operations of combined dense AMLS 
is given in
Table \ref{AMLS: table: overview combined AMLS} where the different tasks of the method are denoted by
$(\ctask 1)$--$(\ctask 8)$. It has to be noted that 
in task (\ctask5)\:b) 
the matrix $Q:=[Q_A,Q_B]$ 
possibly has to be orthogonalised (e.g., by a QR-factorisation) in order to eliminate linear dependent columns in $Q$
which are not wanted when the reduced eigenvalue problem is solved in task (\ctask7).

\begin{figure*}[p]
\centering
\begin{minipage}[c]{0.45\textwidth}
      \begin{overpic}[width=8.2cm]{eigenfunc_1.eps}
      \end{overpic}
\end{minipage}
\hfill
\begin{minipage}[c]{0.45\textwidth}
      \begin{overpic}[width=8.2cm]{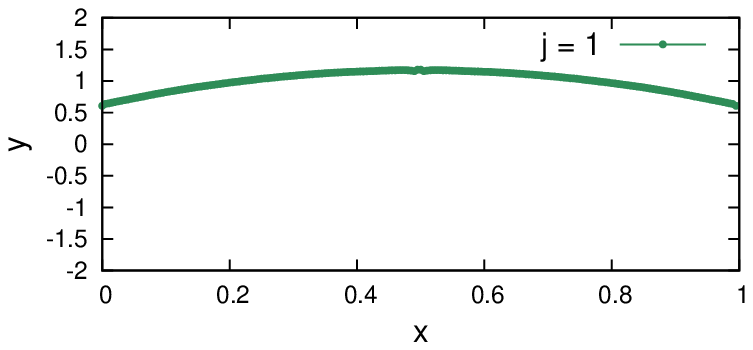}
      \end{overpic}
\end{minipage}
\\
\begin{minipage}[c]{0.45\textwidth}
      \begin{overpic}[width=8.2cm]{eigenfunc_2.eps}
      \end{overpic}
\end{minipage}
\hfill
\begin{minipage}[c]{0.45\textwidth}
      \begin{overpic}[width=8.2cm]{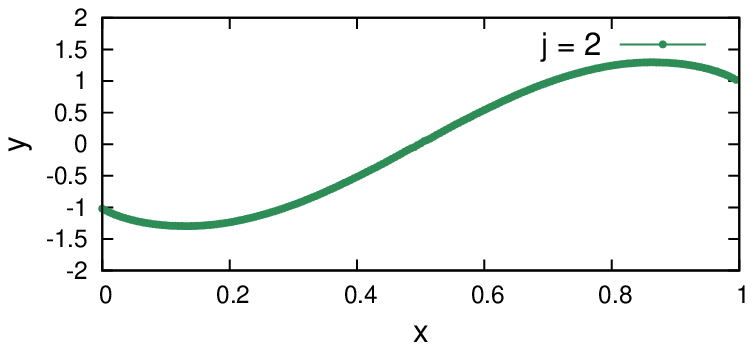}
      \end{overpic}
\end{minipage}
\\
\begin{minipage}[c]{0.45\textwidth}
      \begin{overpic}[width=8.2cm]{eigenfunc_3.eps}
      \end{overpic}
\end{minipage}
\hfill
\begin{minipage}[c]{0.45\textwidth}
      \begin{overpic}[width=8.2cm]{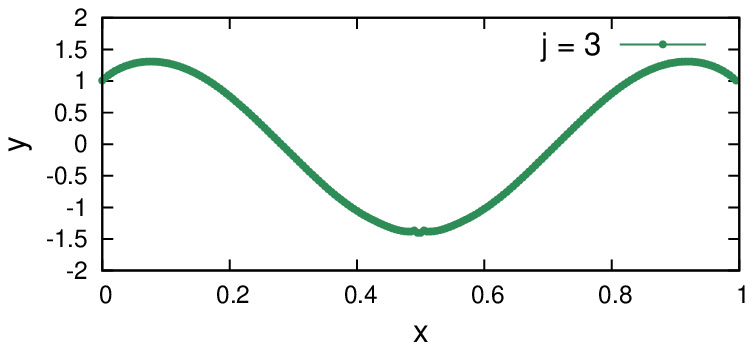}
      \end{overpic}
\end{minipage}
\\
\begin{minipage}[c]{0.45\textwidth}
      \begin{overpic}[width=8.2cm]{eigenfunc_4.eps}
      \end{overpic}
\end{minipage}
\hfill
\begin{minipage}[c]{0.45\textwidth}
      \begin{overpic}[width=8.2cm]{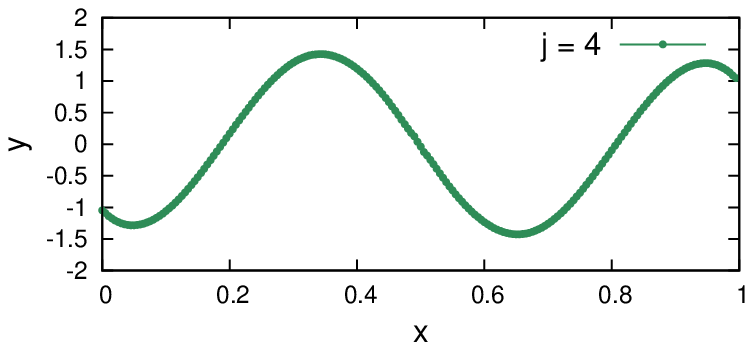}
      \end{overpic}
\end{minipage}
\\
\begin{minipage}[c]{0.45\textwidth}
      \begin{overpic}[width=8.2cm]{eigenfunc_5.eps}
      \end{overpic}
\end{minipage}
\hfill
\begin{minipage}[c]{0.45\textwidth}
      \begin{overpic}[width=8.2cm]{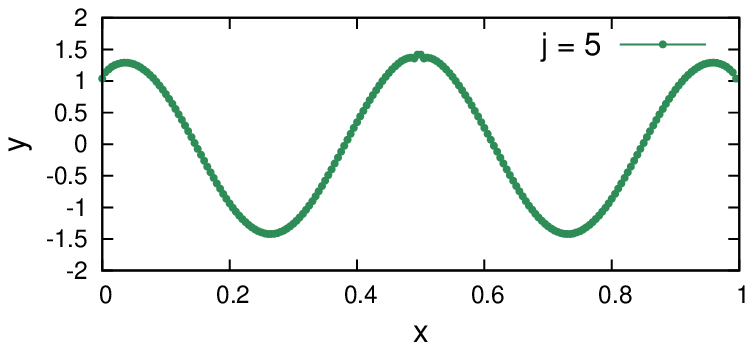}
      \end{overpic}
\end{minipage}
\caption{
Approximations for the eigenfunctions
$u_j$ of problem \eqref{results: eigenvalue problem example}
for $j=1,\ldots,5$:
On the left-hand side are the approximations provided by the discrete 
eigenvectors $x_j\hidx$ of $(K\hidx,M\hidx)$
where the underlying mesh width $h$ leads to a 
discretisation with 200 DOF, and 
on the right-hand side are the approximations that are provided by the Ritz-vectors 
$\widehat{y}_j\hidx$ of the combined dense AMLS method.}
\label{introduction: eigenfunctions classical dense AMLS combi}
\end{figure*}

Using the same parameter setting 
which has been used for the benchmarks of dense AMLS, 
the combined dense AMLS method leads to 
a reduced eigenvalue problem of size $\ksum=10+10=20$.
The resulting accuracy of combined dense AMLS for the eigenvalue approximation 
is displayed in Table
\ref{AMLS: table: combined dense AMLS}.
The benchmarks show that inequality 
\eqref{Results: postulation} is fulfilled for $\nev = 12$, i.e., 
the approximation error of the combined dense AMLS method is of the 
same order as the discretisation error for $\nev = 12$.
Furthermore, Figure \ref{introduction: eigenfunctions classical dense AMLS combi}
shows that the eigenfunction approximations provided by combined dense AMLS are 
nearly identical with the approximations provided by the discrete
eigenvectors $x_j\hidx$ of $(K\hidx,M\hidx)$.

The numerical results show that the combined dense AMLS method is potentially 
well suited for the solution of problem 
\eqref{introduction: discrete EVP}
respectively
\eqref{introduction: continuous EVP classic}.
However, up to this point the combined dense AMLS method is expensive. For example, 
the computation of tasks (\ctask2) and task (\ctask3)
cause costs of the order $\Ocal(N^3)$. To make out of combined dense AMLS an efficient eigensolver 
major improvements have to be made which are described in the following section.

\newpage

\section{Recursive Approach and Use of Hierarchical Matrices}
\label{section: Recursive Approach and Use of Hierarchical Matrices}

There are two fundamental improvements of the combined dense AMLS method which will accelerate the computational time 
of the method dramatically.
These improvements are motivated by the results presented in \cite{Diss_Gerds,Gerds},
which show that the computational performance of classical AMLS 
in the context of elliptic PDE eigenvalue problems can be significantly improved by
using a recursive approach and by using the concept of hierarchical matrices. 
The two fundamental improvements of the combined dense AMLS method are described in the following.

\subsection*{Recursive Approach}
To improve the computational efficiency of the combined dense AMLS method, it is proposed
to apply combined dense AMLS recursively in task (\ctask4).
More precisely, it is proposed to apply combined dense AMLS
recursively
for the solution of the subproblems 
$({K}_{11},{M}_{11})$
and
$(\widetilde{K}_{22},\widetilde{M}_{22})$
until the size $N_i$ ($i=1,2$) of the subproblems is small enough
---
e.g., smaller than 
a given threshold $\MaxDofSubdomain \in \N$
---
to be solved easily by a direct solver.
The depth of the recursion $\rmax \in \N$ depends on the size $N$ of the problem $(K,M)$ and is described by 
\begin{equation}
N \approx 2^{\rmax} \,\MaxDofSubdomain 
\; \Rightarrow \;
\rmax \sim \log N.
\label{hamls: level recursion}
\end{equation}

\subsection*{Usage of Hierarchical Matrices}
The dense matrix operations in task (\ctask2) and task (\ctask3) 
of combined dense AMLS can be performed 
much more efficiently using the so-called hierarchical matrices (short \HM-matrices).
\HM-matrices \cite{First_HMatrix_part_I,Hackbusch_Buch_HMatrix} are matrices which basically consist 
of large low-rank submatrices and few small dense submatrices. This special matrix format
is very well suited to approximate dense matrices which, e.g., arise in the context 
of the discretisation of partial differential and integral equations
\cite{elliptic_H_Inverse,Faustmann_Melenk_Praetorius_H_LU,H_LU}.
The basic idea of \HM-matrices is to reorder the rows and columns of a matrix in a such way that 
certain submatrices can be approximated or represented by low rank matrices.
Using the concept of low rank approximation a fully populated but data-sparse matrix of size 
$N \times N$ can be represented using only $\mathcal{O}(N \log^{\alpha} N)$ data instead of storing $N^2$ entries
where $\alpha \in \{1,\ldots,4\}$ (cf. \cite{Construction_and_arithmetics_hmatrix,DD_H_LU}).
Most importantly, 
the concept of low rank approximation allows the 
\HM-matrices to perform exact matrix-vector multiplication and approximated
matrix(-matrix) operations (e.g., addition, multiplication, inversion, \LU-factorisation) 
in almost linear complexity $O(N \log^{\alpha} N)$.

To use the fast \HM-matrix arithmetic in task (\ctask2) and task (\ctask3) of
combined dense AMLS, the dense stiffness matrix 
$K$ and the sparse mass matrix $M$ have first to be 
approximated by \HM-matrices. How this is done is briefly described in the following for the matrix $K$.
The same approach has to be used for the matrix $M$. To see more details
it is referred, e.g., to \cite{Gerds,Construction_and_arithmetics_hmatrix}.
\\

Since  matrix $K \in \R^{N\times N}$ results from a finite element discretisation
each row and column index $i \in I:= \{ 1, \ldots, N \}$ is associated
with a basis function $\varphi_i\hidx$ of the underlying finite element space $V_h$. The support of each index set 
$t \subset I$ is defined by
\begin{equation*}
\Omega_t:= 
\bigcup_{i \in t} 
\;
\mathrm{supp} \bigl ( \varphi_i\hidx \bigr ),
\end{equation*}
and correspondingly each submatrix
\begin{equation*}
K|_{s \times t} 
:= 
\left (
K_{ij}
\right )_{i \in s, j \in t}
\qquad
\textnormal{with } 
s,t \subset I
\end{equation*}
of $K$ is associated with geometry information.
Based on the geometric separation of the supports $\Omega_s$ and $\Omega_t$ certain subblocks 
$s \times t \subset I \times I$ can be identified that allow
a low rank approximation of the respective submatrices $K|_{s \times t}$.
More precisely, submatrices $K|_{s \times t}$ whose index sets $s$ and $t$ fulfil the 
so-called \emph{admissibility condition}
\begin{equation}
\min
\bigl \{
\mathrm{diam}( \Omega_s ),\mathrm{diam}( \Omega_t )
\bigr \}
\leq
\eta \; \mathrm{dist} ( \Omega_s,\Omega_t )
\label{HMatrix: admissibility condition}
\end{equation}
are well suited for a low rank approximation (cf. \cite{Construction_and_arithmetics_hmatrix}).
The quantities
\begin{alignat*}{2}
\mathrm{diam}( \Omega_t )
&:=
\max
& & \bigl \{
\| x-y \|_2
\; : \;
x, y \in \Omega_t
\bigr \},
\rule[-3mm]{0mm}{3mm}
\\
\mathrm{dist} ( \Omega_s,\Omega_t )
&:=
\min 
& & \bigl \{ 
\| x-y \|_2
\; : \;
x \in \Omega_s, y \in \Omega_t
\bigr \}
\end{alignat*}
are the diameter and the distance of the supports $\Omega_s$ and $\Omega_t$, and 
the parameter $\eta > 0$ controls the number of admissible subblocks $s \times t$
and is typically set to $\eta = 1$ (see, e.g., \cite{Construction_and_arithmetics_hmatrix}).

Subblocks $s \times t$ of $I \times I$ which fulfil the admissibility condition 
\eqref{HMatrix: admissibility condition} are
called in the following \emph{admissible} and the corresponding submatrices $K|_{s\times t}$ are 
approximated by the so-called \Rk-matrices.
\begin{definition}[\Rk-matrix]
Let $k,m,n \in \N_0$. A matrix $R \in \R^{n \times m}$ is called \Rk-matrix if it is factorised by
\begin{equation}
R = A B^T
\; \;
\textnormal{ with suitable matrices }
\; A \in \R^{n \times k}
\textnormal{ and }
B \in \R^{m \times k}.
\label{HMatrix: Rk-matrix}
\end{equation}
\end{definition}
The representation of an \Rk-matrix $R \in \R^{n \times m}$
in factorised form \eqref{HMatrix: Rk-matrix}
is much cheaper than in full-matrix representation
when the rank $k$ is small compared to $n$ and $m$ because only $k(n+m)$ entries have to be
stored instead of $n \, m$. Furthermore, when $k$ is small the product and the sum of \Rk-matrices can be evaluated 
much more efficiently than in full-matrix representation.

To exploit the low rank approximation property of submatrices $K|_{s \times t}$ fulfilling 
\eqref{HMatrix: admissibility condition} the row and column indices of $K$ have to be reordered.
For this purpose the index set $I$ is divided according to a geometric bisection of its support $\Omega_I$ into two
disjoint index sets $s,t \subset I$ with $I = s \, \dot{\cup} \, t$. In this context we denote $s$ and $t$ 
as the \emph{sons} of $I$ and $S(I):= \{ s, t \}$ as the \emph{set of sons} of $I$. 
This geometric bisection of $I$ is applied recursively to the son index sets until 
the cardinality of an resulting index set falls below
some given threshold $n_{min} \in \N$. This recursive geometric bisection of the index set $I$ results in a 
disjoint partition of $I$ where the obtained subsets of the partitioning tend 
to be geometrically separated
(see, e.g., \cite{Construction_and_arithmetics_hmatrix} for details).
Finally, the row and column indices of the matrix $K$ are reordered correspondingly to the performed 
partitioning of $I$.

Given the partitioning of $I$ and 
the admissibility condition \eqref{HMatrix: admissibility condition} 
the \HM-matrix approximation of $K$ is computed by applying Algorithm 
\ref{HMatrix: alg: H-Matrix Approximation} to the block index set $I \times I$ and to the reordered matrix $K$. 
The resulting \HM-matrix approximation is denoted by $K^{\HM}$ in the following.
Using Algorithm \ref{HMatrix: alg: H-Matrix Approximation}
the block index set $I \times I$ is recursively
subdivided into subblocks $s \times t$ until the subblock gets admissible or the size of the subblock falls
below the threshold $n_{min}$.
Submatrices $K|_{s \times t}$ of admissible blocks $s \times t$ are approximated by \Rk-matrices
and submatrices of inadmissible blocks are represented in the full-matrix format.
To control 
in Algorithm \ref{HMatrix: alg: H-Matrix Approximation}
the approximation quality of the \Rk-matrix approximation the fixed rank $k$ is replaced by
an adaptive rank: 
Each submatrix $K|_{s\times t}$
corresponding to an admissible subblock $s \times t$ can be approximated by an \Rk-matrix $R$ such 
that 
\begin{equation}
\frac{\| \, K|_{s\times t} - R \; \|_2}
{\| \, K|_{s\times t} \; \|_2}  \leq \varepsilon
\label{HMatrix: approximation accuracy}
\end{equation}
where $\varepsilon > 0$ is some arbitrary prescribed approximation accuracy
and where the rank $k \in \N_0$ is as small as possible (cf. \cite{Construction_and_arithmetics_hmatrix}).

\begin{algorithm}[h]
\caption{Computation of the \HM-Matrix Approximation}
\label{HMatrix: alg: H-Matrix Approximation}

\begin{small}
\begin{algorithmic}
\Procedure{GetHMatrixApproximation}{$K$, $\varepsilon$, $n_{min}$, $s \times t$}
  \If{$s \times t$ is admissible}
    \State{$K^{\HM}|_{s \times t}:=$ \Rk-matrix approximation of $K|_{s \times t}$ with accuracy $\varepsilon$;}
  \ElsIf{$\min \{ \# s, \# t \} \leq n_{min}$} 
    \State{$K^{\HM}|_{s \times t}:=$ full-matrix representation of $K|_{s \times t}$;} 
    \Comment{$n_{min}$ affects the minimal size of submatrices} 
  \Else
    \State{$S(s \times t) := \bigl \{ s' \times t' \; | \; s' \in S(s),\, t' \in S(t) \bigr \}$;}
    \Comment{$S(t)$ denotes the set of sons of $t \subset I$}
    \ForAll{$s' \times t' \in S(s \times t)$}
      \State{$K^{\HM}|_{s' \times t'}:=$\Call{GetHMatrixApproximation}{ $K$, $\varepsilon$,  $n_{min}$, $s' \times t'$ };}
    \EndFor
  \EndIf
\EndProcedure
\end{algorithmic}
\end{small}

\end{algorithm}

Using the same approach, we compute the \HM-matrix approximation $M^{\HM}$ of $M$, i.e., 
first we reorder the row and column indices of $M$ according to the partitioning of $I$,
and then we apply Algorithm \ref{HMatrix: alg: H-Matrix Approximation} to 
the block index set $I \times I$ and 
the reordered matrix $M$.
The approximation errors ${\| {K} - {K}^{\HM} \|_2}$ and ${\| {M} - {M}^{\HM} \|_2}$ are controlled by the 
accuracy $\varepsilon$ of the \Rk-matrix approximation in \eqref{HMatrix: approximation accuracy}.

The \HM-matrix approximations $K^{\HM}$ and $M^{\HM}$ of $K$ and $M$ can now be used to compute
in task (\ctask2) and (\ctask3) of combined dense AMLS 
the transformed eigenvalue problem $(\widetilde{K},\widetilde{M})$. 
Note that the \HM-matrices $K^{\HM}$ and $M^{\HM}$ possess a recursive $2\times2$ block partitioning of the form
\begin{equation}
{K}^{\HM}
= 
\begin{bmatrix}
{K}^{\HM\rule[-0.5mm]{0mm}{0mm}}_{11}&    {K}^{\HM\rule[-0.5mm]{0mm}{0mm}}_{12} \rule[-3mm]{0mm}{0mm}\\
{K}^{\HM\rule[-0.5mm]{0mm}{0mm}}_{21}&    {K}^{\HM\rule[-0.5mm]{0mm}{0mm}}_{22}
\end{bmatrix}
\qquad \textnormal{and} \qquad
{M}^{\HM}
= 
\begin{bmatrix}
{M}^{\HM\rule[-0.5mm]{0mm}{0mm}}_{11}&    {M}^{\HM\rule[-0.5mm]{0mm}{0mm}}_{12} \rule[-3mm]{0mm}{0mm}\\
{M}^{\HM\rule[-0.5mm]{0mm}{0mm}}_{21}&    {M}^{\HM\rule[-0.5mm]{0mm}{0mm}}_{22}
\end{bmatrix}
\label{Hmatrix: marix partitioning}
\end{equation}
where ${K}^{\HM\rule[-0.5mm]{0mm}{0mm}}_{ij},{M}^{\HM\rule[-0.5mm]{0mm}{0mm}}_{ij} \in \R^{N_i \times N_j}$
are either full-matrices, \Rk-matrices or \HM-matrices which again have a $2\times 2$ block partitioning.
Because of the matrix partitioning \eqref{Hmatrix: marix partitioning} of $K^{\HM}$ and $M^{\HM}$
(which is based on a geometric bisection of the support $\Omega_I$), 
there is no need to perform the domain substructuring and the associated matrix 
partitioning of task (\ctask1) in combined dense AMLS.
Using the fast \HM-matrix arithmetic and the \HM-matrix approximations $K^{\HM}$ and $M^{\HM}$,
the block diagonalisation of $K$ and the transformation of $M$ can be computed 
very efficiently via
\begin{equation}
K^{\HM} \approx L^{\HM} \widetilde{K}^{\HM} (L^{\HM})^T
\quad
\textnormal{ and }
\quad
\widetilde{M}^{\HM}
\hspace*{-1mm}
\approx
({L}^{\HM})^{-1} 
M^{\HM}
({L}^{\HM})^{-T}
\label{HAMLS: transformation K and M}
\end{equation}
in $\mathcal{O}(N\log^{\alpha} N)$.
The corresponding algorithm for the computation of 
\eqref{HAMLS: transformation K and M}
is based on a recursive approach which is applied block-wise to the matrix structure
\eqref{Hmatrix: marix partitioning}, which exploits the \Rk-matrix 
representation of submatrices fulfilling \eqref{HMatrix: admissibility condition}, and 
which applies the inexpensive addition and multiplication of \Rk-matrices
(see, e.g., \cite{Construction_and_arithmetics_hmatrix} for details).
The \HM-matrix operations in \eqref{HAMLS: transformation K and M} are performed not exactly
but only approximatively, however, the approximation errors
${\| \widetilde{K} - \widetilde{K}^{\HM} \|_2}$ and 
${\| \widetilde{M} - \widetilde{M}^{\HM} \|_2}$
can be controlled by the chosen accuracy $\varepsilon$ in \eqref{HMatrix: approximation accuracy}.

\vspace*{10mm}


The improved version of combined dense AMLS with the recursive approach and the use of the fast \HM-matrix arithmetic 
is simply called \emph{dense \HM-AMLS} in the following.
It is proposed to use the \HM-matrix software library HLIBpro
\cite{Intro_HMatrix,HLIBpro,Parallel_HMatrix_shared_memory} for the implementation of the method;
and since dense \HM-AMLS has much in common with 
the method described in \cite{Diss_Gerds}, it is recommended to see \cite[Section 8]{Diss_Gerds} for corresponding 
implementation issues.

Beside the number of sought eigenpairs $\nev$ and the number of degrees of freedom $N$ 
the computational costs of dense \HM-AMLS 
depend on the chosen accuracy $\varepsilon$ of the approximative \HM-matrix operations in 
\eqref{HAMLS: transformation K and M},
and on the applied modal truncation in task (\ctask4), i.e, the number of selected 
eigenvectors ${k}_i$ in \eqref{AMLS classic: partial eigensolution}.
Coarsening the accuracy $\varepsilon$ of the \HM-matrix operations and decreasing 
the number $k_i$ of selected eigenvectors 
result in faster computations and reduced memory requirements of dense \HM-AMLS.
On the other side the accuracy $\varepsilon$ has to be fine enough and the size $k_i$ ($i=1,2$) 
has to be large enough so that the error caused by the approximative \HM-matrix arithmetic and the error caused
by the modal truncation are of the same order as the discretisation error. If this is the case the approximation
quality of dense \HM-AMLS competes with the approximation quality of a classical approach
[cf. Remark \ref{AMLS classic: remark: mode selection and reduced problem}\:i)].
Hence, the aim arises to select the parameters $\varepsilon$ and $k_i$ in such a way that dense \HM-AMLS
reaches the accuracy of a classical approach while the computational costs of \HM-AMLS are reduced as much 
as possible.

In the following the computational costs of dense \HM-AMLS are discussed in detail.
The discussion is restricted to problems
where $\nev \sim N^\beta$ eigenpairs are sought with some $\beta \in (0,1/3]$.
In task (\ctask4) the mode selection strategy proposed 
by Remark \ref{AMLS: remark: mode selection} is applied, i.e., 
the eigenvectors of 
$({K}_{11}^{\HM},{M}_{11}^{\HM})$
and
$(\widetilde{K}_{22}^{\HM},\widetilde{M}_{22}^{\HM})$ associated to the $k_i$ (i=1,2)
smallest eigenvalues are computed where
$k_i \in \Ocal(N_i^{\beta})$, and thus it is guaranteed that the size $\ksum$ of the 
reduced problem in task (\ctask7) is bounded by $\Ocal(\nev)$.
Using this setting the costs of dense \HM-AMLS can be summarised as follows:
\begin{itemize}

\item 
The computational costs for the matrix partitioning in (\ctask1), respectively 
the costs for the computation of the \HM-matrix approximations
$K^{\HM}$ and $M^{\HM}$ are of the order ${\cal O}(N \log N)$.
The computational costs for task $(\ctask 2)$ and $(\ctask 3)$ 
using the fast \HM-matrix arithmetic
are of the order
$\mathcal{O}(N\log^{\alpha} N)$.

\item
In task (\ctask4) four subproblems arise 
(cf. Table \ref{AMLS: table: overview AMLS})
of size $N_i \approx N/2$ (with $i=1,2$).
In the case that the subproblems are small enough 
(i.e., it holds $N_i \leq \MaxDofSubdomain$) 
the subproblems are solved by a direct solver 
leading to costs of the order $\Ocal(1)$, otherwise the subproblems 
are solved recursively by the dense \HM-AMLS method.

\item
The computation of the matrices $Q_A=(L^{\HM})^{-T}Z$ and 
$Q_B=(L^{\HM})^{-T}Z$ in task (\ctask5) can be performed efficiently 
by exploiting the fast \HM-matrix times vector multiplication
[backward substitution in $(L^{\HM})^T$]. Since the number of columns of $Q_A \in R^{N\times \ksum_A}$
and $Q_B \in R^{N\times \ksum_B}$ is bounded by $\ksum_A,\ksum_B \in \Ocal( N^{\beta})$, 
the computation of $Q_A$ and $Q_B$ can be performed in 
${\cal O} (  N^{\beta} \, N \log^{\alpha} N )$.
The computation of the matrix $Q \in \R^{N\times\ksum}$ where $\ksum \leq \ksum_A+\ksum_B$,
more precisely, the elimination 
of linear dependent columns in the matrix $[Q_A,Q_B]$ by orthogonalisation, leads to computational costs of the order
$\Ocal (N^{2\beta}N)$. 
Since $\nev \sim N^{\beta} $ the
total costs of task (\ctask5) can be bounded by
\begin{equation*}
{\cal O} \bigl ( \, N^{\beta} \, N\log^{\alpha} N \, \bigr ) 
+ \,
{\cal O} \bigl (  \, N^{2\beta} \, N   \, \bigr )
= 
{\cal O} \bigl ( \, \nev \, N\log^{\alpha} N \, \bigr ) 
+ \,
{\cal O} \bigl (  \, \nev^2 \, N   \, \bigr ).
\end{equation*}

\item
The \HM-matrix structure of ${K}^{\HM}$ and ${M}^{\HM}$ can be exploited as well in task (\ctask6)
using the fast \HM-matrix-vector multiplication for the computation of the reduced matrices 
${\widehat{K}}, {\widehat{M}} \in \R^{\ksum \times \ksum}$: 
The multiplications 
$Q^T({K}^{\HM}Q)$ and 
$Q^T({M}^{\HM}Q)$ involve in total
$2\ksum$ \HM-matrix times vector  multiplications with costs of the order 
${\cal O}(N \log^{\alpha}N)$ 
plus $2\ksum^2$ scalar products of length $N$.
Hence, the costs of (\ctask6) sum up to 
\begin{equation*}
{\cal O} \bigl ( \, \ksum \, N \log^{\alpha} N \,+\, \ksum^2 \,N \, \bigr ) 
\; \leq \;
{\cal O} \bigl ( \, \nev \, N \log^{\alpha} N \,+\, \nev^2 \,N   \, \bigr ).
\end{equation*}

\item
Since the size $\ksum$ of the reduced problem 
$(\widehat{K},\widehat{M})$ is bounded by 
$\ksum \in \Ocal(\nev)$ and since
we aim at $\nev \in \Ocal(N^{1/3})$ eigenpairs, 
we can use in task (\ctask7) a dense linear algebra solver with cubic complexity for the solution of the reduced 
problem whereby the computational costs still remain in ${\cal O}(N)$.

\item
Finally, in task $(\ctask 8)$ in total $N\nev$ scalar products of length $\ksum$
have to performed leading to costs of the order
\begin{equation*}
{\cal O} \bigl ( \, \nev \ksum N \, \bigr ) 
\; \leq \;
{\cal O} \bigl ( \, \nev^2 N \, \bigr ).
\end{equation*}
\end{itemize}
Hence, the computational costs of dense \HM-AMLS ---
without the costs of
the recursive calls of dense \HM-AMLS  in task (\ctask4)
---
are bounded by
${\cal O} ( \nev N \log^{\alpha}N + \nev^2 N )$.
Note that when dense \HM-AMLS is applied recursively to a problem of size $N$, 
that then on the next level of the 
recursion four new subproblems
(cf. Table \ref{AMLS: table: overview AMLS})
arise of approximate size $N/2$ which are handled recursively
by dense \HM-AMLS.
It follows that the overall computational costs of dense \HM-AMLS, including the costs for the 
recursive calls, are bounded by
\begin{align*}
&
{\cal O}\bigl ( \, \nev \, N \log^{\alpha}N + \nev^2 \, N \, \bigr )
+
\sum_{l=1}^{\rmax}
4^l \;
\Ocal \biggl ( 
\, \Bigl ( \dfrac{N}{2^l} \Bigr )^{\beta} \, 
\dfrac{N}{2^l} \log^{\alpha}\Bigl ( \dfrac{N}{2^l} \Bigr ) 
\, + \, 
\Bigl ( \dfrac{N}{2^l} \Bigr )^{2\beta} \, \dfrac{N}{2^l} \, 
\biggr )  
\\
\; = \;
&
{\cal O}\bigl ( \, \nev \, N \log^{\alpha}N + \nev^2 \, N \, \bigr )
+
\sum_{l=1}^{\rmax}
2^l \;
\Ocal \biggl ( 
\, \Bigl ( \dfrac{N}{2^l} \Bigr )^{\beta} \, 
N \log^{\alpha} N
\, + \, 
\Bigl ( \dfrac{N}{2^l} \Bigr )^{2\beta} \, N \, 
\biggr ) 
\\
\; = \;
&
{\cal O}\bigl ( \, \nev \, N \log^{\alpha}N + \nev^2 \, N \, \bigr )
+
\sum_{l=1}^{\rmax}
2^l \;
\Ocal \biggl ( 
\, \dfrac{\nev}{2^{l\beta}} \, 
N \log^{\alpha} N
\, + \, 
\dfrac{\nev^2}{2^{2l\beta}} \, N \, 
\biggr ) 
\\
\; = \;
&
{\cal O}\bigl ( \, \nev \, N \log^{\alpha}N \, \bigr )
\sum_{l=0}^{\rmax}
 2^{l(1-\beta)}
\; + \;
{\cal O}\bigl ( \, \nev^2 \, N \, \bigr )
\sum_{l=0}^{\rmax}
2^{l(1-2\beta)}
\\
\; = \;
&
{\cal O}\bigl ( \, \nev \, N \log^{\alpha}N \, \bigr )
\,
\dfrac{
2^{(\rmax+1)(1-\beta)}-1
}{2^{1-\beta}-1}
\; + \;
{\cal O}\bigl ( \, \nev^2 \, N \, \bigr )
\,
\dfrac{
2^{(\rmax+1)(1-2\beta)}-1
}{2^{1-2\beta}-1}
\\
\; \mathop{=}^{\eqref{hamls: level recursion}} \;
&
{\cal O}\bigl ( \, \nev \, N \log^{\alpha}N \, N^{(1-\beta)} \, \bigr )
\; +\;
{\cal O}\bigl ( \, \nev^2 \, N \, N^{(1-2\beta)} \, \bigr )
\;=\;
{\cal O} \bigl ( \, N^2 \log^{\alpha}N + N^2 \, \bigr )
\end{align*}

We can sum up that the theoretical complexity
of dense \HM-AMLS is bounded%
\footnote{Note that the upper bound 
${\cal O} ( N^2 \log^{\alpha}N + N^2 )$
for the computational costs of dense \HM-AMLS 
computing the $\nev \sim N^{\beta}$ smallest eigenpairs
is obtained for all $\beta \in (0,1)$, i.e., also in the case 
when it holds $\beta > 1/3$.
}
by ${\cal O} ( N^2 \log^{\alpha}N + N^2 )$.
The costs of dense \HM-AMLS are theoretically dominated by the costs of 
of the scalar products occurring in tasks (\ctask5), (\ctask6) and (\ctask8)
which are accumulating to $\Ocal(N^2)$, and the costs of the 
\HM-matrix times vector multiplications in tasks (\ctask5) and (\ctask6) which 
are accumulating to ${\cal O} ( N^2 \log^{\alpha}N )$.
However, since scalar products can be computed with peak performance
on today's workstations and compute servers, 
it is expected that these operations will not dominate the overall costs of dense \HM-AMLS in practice.
It is expected as well that in practice the 
\HM-matrix times vector multiplications are as well harmless since the logarithms and constants involved in 
the \HM-matrix vector multiplications are much smaller than for the \HM-matrix operations performed in 
task (\ctask2) and task (\ctask3). Overall, it is expected that in practice the costs 
of dense \HM-AMLS are dominated by the costs of the \HM-matrix operations in task (\ctask2) and (\ctask3) 
which, however,
accumulate as well, due to the recursive calls, to costs of the order ${\cal O} ( N^2 \log^{\alpha}N )$. 
But possibly these costs can be bounded by some approach which limits the number of needed problem 
transformations [task (\ctask2) and (\ctask3)] in the recursive calls. In \cite[Section 8]{Diss_Gerds}
an approach is described for a quite similar method which allows to replace all needed 
problem transformations of the recursive calls by only one global problem transformation.
Furthermore, the very high efficiency of the method described in \cite{Diss_Gerds,Gerds}
(which combines in the context of elliptic PDE eigenvalue problems the classical AMLS method with 
a recursive approach and the usage of the \HM-matrices) 
towards classical approaches motivates that the very similar dense \HM-AMLS method reaches in practice 
possibly the same superior efficiency.

\section{Conclusion}
\label{section: Conclusion}

While domain decomposition methods, like the very efficient AMLS method, are available for the solution 
of elliptic PDE eigenvalue problems, domain decomposition techniques are not known for the solution
integral equation eigenvalue problems. 
To the best of the author's knowledge this paper introduces the very first domain decomposition method for 
integral equation eigenvalue problems. The new method, which is motivated by the classical 
AMLS method, decomposes a global problem with $N$ degrees of freedom into four subproblems each with approximately $N/2$
degrees of freedom. The eigensolutions of the four subproblems are then used to form a subspace
which approximates the sought eigensolutions of the global problem. This domain decomposition technique is called combined dense AMLS method 
and shows very promising results concerning the approximation quality. 
To improve the computational efficiency of the combined dense AMLS method,  
it is proposed to use a recursive approach and the fast \HM-matrix arithmetic for the needed problem transformation.
The theoretical computational complexity of this approach is of the order ${\cal O} ( N^2 \log^{\alpha}N )$, however,
motivated by the results for a quite similar method (cf. \cite{Diss_Gerds,Gerds})
it is expected that the method reaches in practice a much better computational efficiency.
Furthermore, the promising results of the combined dense AMLS method 
show that domain decomposition techniques are also applicable to integral equation eigenvalue problems
and motivate to refine the domain decomposition techniques of combined dense AMLS. 
Possibly the combined dense AMLS method can be refined in such a way that the global problem is 
decomposed into only two subproblems instead of four which might decrease the computational complexity.

\bibliographystyle{abbrv} 
\bibliography{references}

\end{document}